\documentclass[12pt]{article}
\usepackage[utf8]{inputenc} 
\usepackage[english,russian]{babel}
\usepackage{amsmath, amssymb, latexsym}
\usepackage{tikz}
\usepackage{ifthen}
\usepackage{hyperref}
\title{MDS codes in Doob graphs}
\author{Evgeny A. Bespalov, Denis S. Krotov}
\date{}

\makeatletter
\def\@seccntformat#1{\csname the#1\endcsname.\ } 
\makeatother

\newif\ifNoRemark
\def\addtheorem#1#2#3#4{
\ifthenelse{\equal{#2}{}}{}%
{\ifthenelse{\expandafter\isundefined\csname the#2\endcsname}{\newcounter{#2}}{}}
\newenvironment{#1}[1][\global\NoRemarktrue]
{\par\addvspace{2mm plus 0.5mm minus 0.2mm}\noindent 
{\bf #3}\ifthenelse{\equal{#2}{}}{}%
{\refstepcounter{#2}{\bf ~\csname the#2\endcsname}}%
{\bf \vphantom{##1}\ifNoRemark.\ \else\ (##1).\fi}\begingroup #4}%
{\endgroup\par\addvspace{1mm plus 0.5mm minus 0.2mm}\global\NoRemarkfalse}
\expandafter\newcommand\csname b#1\endcsname{\begin{#1}}
\expandafter\newcommand\csname e#1\endcsname{\end{#1}}
}

\addtheorem{theoreman}{thrm}{Теорема}{\sl}
\addtheorem{lemman}{lmm}{Лемма}{\sl}
\addtheorem{predln}{prpstn}{Педложение}{\sl}
\addtheorem{coroll}{crllr}{Следствие}{\sl}
\addtheorem{remarkn}{rmrk}{Замечание}{}
\addtheorem{definitionn}{def}{Определение}{}

 \newenvironment{proof}[1][\hspace{-1.0ex}]%
  {\par\addvspace{1mm}{\sc Доказательство\hspace{1.0ex}{#1}.} }%
  {\quad$\blacktriangle$\par\addvspace{1mm}}
	
\def\shpart#1 #2 #3 #4!{node [#1] {} +(1,0) node [#2] {} +(2,0) node [#3] {} +(3,0) node [#4] {}}
\def\sh#1#2#3#4{\begin{tikzpicture}[
scale=0.7,
nz/.style={circle,fill=white,draw=black, 
           inner sep=2.5pt},
xz/.style={circle,fill=black!50!white,draw=black, 
           inner sep=2.5pt},
vz/.style={circle,fill=black!23!white,draw=black, 
           inner sep=2.5pt},
wz/.style={circle,fill=black!66!white,draw=black, 
           inner sep=2.5pt},
zz/.style={circle,fill=black,draw=black, 
           inner sep=2.5pt},
scale=0.7]
\begin{scope}
\clip [xslant=-0.577] (-1.4,-1.20) rectangle (2.4,2.1);
\draw[xslant=0.577,ystep=.866,xstep=1,draw=black] (-4.9,-2.1) grid (5.4,3.9);
\draw[xslant=-0.577,ystep=9.866,xstep=1,draw=black] (-3.4,-2.1) grid (6.4,3.9);
\draw (-120:1) \shpart #4!
++(120:1) \shpart #3!
++(120:1) \shpart #2!
++(120:1) \shpart #1!;
\end{scope}
\end{tikzpicture}}

\begin{document}
\title{МДР коды в графах Дуба%
\thanks{Исследование выполнено за счет гранта Российского научного фонда (проект №14-11-00555).
Адрес авторов: Институт математики им. С. Л. Соболева СО РАН, проспект Академика Коптюга
4, Новосибирск, 630090, Россия (e-mail: bespalovpes@mail.ru, krotov@math.nsc.ru).}
}
\author{Беспалов Е. А., Кротов Д. С.}

\maketitle
\begin{abstract}
The Doob graph $D(m,n)$, where $m>0$, is the direct product of $m$ copies of The Shrikhande graph and $n$ copies of the complete graph $K_4$ on $4$ vertices. 
The Doob graph $D(m,n)$ is a distance-regular graph with the same parameters as the Hamming graph $H(2m+n,4)$. 
In this paper we consider MDS codes in Doob graphs with code distance $d \ge 3$. 
We prove that if $2m+n>6$ and $2<d<2m+n$, then there are no MDS codes with code distance $d$. 
We characterize all MDS codes with code distance $d \ge 3$ in Doob graphs $D(m,n)$ when $2m+n \le 6$.
We characterize all MDS codes in $D(m,n)$ with code distance $d=2m+n$ for all values of $m$ and $n$.

Граф Дуба $D(m,n)$, где $m>0$, декартово произведение $m$ копий графа Шрикханде и $n$ копий полного графа $K_4$ на 4 вершинах. 
Граф Дуба $D(m,n)$ дистанционно-регулярный граф с теми же параметрами, что и граф Хэмминга $H(2m+n,4)$. В работе рассматриваются МДР коды в графах Дуба $D(m,n)$ с кодовым расстоянием $d \ge 3$. Доказано, что если $2<d<2m+n$ и $2m+n>6$, то в графе $D(m,n)$ не существует МДР кодов с кодовым расстоянием $d$.
Характеризованы все МДР коды с кодовым расстоянием $d \ge 3$ в графах $D(m,n)$ при $2m+n \le 6$. 
Также характеризованы все МДР коды в графе $D(m,n)$ с кодовым расстоянием $d=2m+n$.
 \end{abstract}

\def\VV{{\scriptscriptstyle\mathrm{V}}}


\section{Введение}\label{s:intro}

В этой работе мы изучаем МДР коды в графах Дуба с кодовым расстоянием $d \ge 3$.

$(N, 4^k,d)$ МДР кодом в графе Хэмминга $H(N,q)$ называется код мощности $4^k$ с кодовым расстоянием $d=n-k+1$. Исследование МДР кодов в графах Хэмминга является важной областью в теории кодирования. Подробнее об этом можно посмотреть в \cite{MWS}. В общем случае вопрос существования и классификации МДР кодов с данными параметрами остается открытой проблемой. О классификации МДР кодов при небольших зачениях $q$ можно посмотерть например в \cite{KokOst:Gr-Lat} и \cite{KKO:smallMDS}.

Случай $q=4$ особый для графов Хэмминга $H(N,q)$, 
так как только в этом случае граф $H(N,4)$, $N \ge 2$, не определяется как дистанционно-регулярный граф с данными параметрами.
Другой граф с теми же параметрами, что и граф Хэмминга $H(2m+n,4)$, это граф Дуба $D(m,n)$, $N=2m+n$.
МДР коды в графах Хэмминга $H(N,4)$ изучаются в работах \cite{Alderson:MDS4} и \cite{KroPot:4}.

$(m+n,4^k,d)$ МДР кодом называется код в графе $D(m,n)$ с кодовым расстоянием $d$ и мощностью кода $4^k$.
В этой статье характеризованы все МДР коды с кодовым расстоянием $d \ge 3$. 
Если в графах $H(N,4)$ при $k=1$ существует единственный с точностью до эквивалентности $(n,4^1,d)$ МДР код, то в графе $D(m,n)$ таких кодов будет значительно больше.   
При $2m+n \le 6$ найдено число МДР кодов с кодовым расстоянием $d=3,4$, и все такие коды с точностью до эквивалентности приведены в Приложении.
Доказано, что при $2m+n>6$ и $2<d<2m+n$ не сущестует $(m+n,4^k,d)$ МДР кодов. Стоит заметить, что для МДР кодов в графах Хэмминга $H(n,q)$ существует предположение, что линейные МДР коды существует только при $n \le q+2$, т.е. при $q=4$ длина кода не превосходит 6. Подробнее об этом смотри в \cite{Ball:2012:1} и \cite{BalBeu:2012}. В полученных результатах для МДР кодов в графах Дуба $D(m,n)$ значение $2m+n$ также не превосходит 6.   

В разделе \ref{s:aux} мы приводим необходимые определения и вспомогательные утверждения.
В разделе \ref{s:main} в теореме \ref{t:main} сформулирован основной результат, в котором приведены значения для числа классов эквивалентности МДР кодов с кодовым расстоянием $d \ge 3$. Доказательство этой теоремы разбито на несколько частей, сформулированных в виде предложений \ref{p:mn1}-\ref{p:notexist}, в зависимости от параметров кода. 
В разделе \ref{s:mn1} рассматриваются коды с параметром $k=1$.  
В разделе \ref{s:mn2} рассматриваются коды с параметрами $k=2$ и $d=3,4$, а именно с параметрами 
$(2+0,4^2,3)$, $(1+2,4^2,3)$, $(2+1,4^2,4)$ и $(1+3,4^2,4)$.
В разделе \ref{s:mn33} рассматриваются коды с параметрами $k=3$, $d=3,4$.
В разделе \ref{s:notexist} доказывается, что при $2 < d < 2m+n$ и $2m+n>6$ не существует $(m+n,4^k,d)$ МДР кодов.
В приложении приведены все МДР коды при $4 \le 2m+n \le 6$ и $d=3,4$ с точностью до эквивалентности. 

\section{Определения и вспомогательные утверждения}\label{s:aux}
\emph{Граф Шрикханде} $\mathrm{Sh}$ --- это граф Кэли над группой $Z^2_4$ с порождающим множеством $\{01,03,10,30,11,33\}$ (вершины графа --- элементы группы $Z^2_4$, которые мы обозначим $00, 01, 02, \ldots, 33$; две 
вершины смежны тогда и только тогда, когда их разность принадлежит пораждающему множеству).
Полный граф $K=K_4$ --- граф Кэли над группой $Z_4$  с порождающим множеством $\{1,2,3\}$. 
 

Пусть $m,n$ неотрицательные целые числа.  Обазначим через $D(m,n)=Sh^m \times K^n$ граф, являющийся прямым произведением $m$ копий графа Шрикханде и $n$ копий 
полного графа $K_4$. 
Если $m>0$, то такой граф называется \emph{графом Дуба}, тогда как $D(0,n)$ --- граф Хэмминга $H(n,4)$. Граф Дуба --- дистанционно-регулярный с теми же параметрами, что и граф Хэмминга $H(2m+n,4)$. 

Множество вершин графа $G$ обозначим через $\VV G$.
Вершины графа $D(m,n)$ мы будем обозначать через $(s_1,\ldots,s_m;h_1,\ldots,h_n)$, 
где $s_i\in Z^2_4$,
$h_i \in Z_4$.
Для произвольных вершин $a,b \in \VV D(m,n)$ обозначим через $d(a,b)$ расстояние между вершинами $a,b$ в графе $D(m,n)$. 
Так как $D(m,n)$ декартово произведение графов $Sh$ и $K$, то для двух вершин $c=(s_1,\ldots,$ $s_m;h_1,\ldots,h_n)$ и $c'=(s'_1,\ldots,s'_m,h_1,\ldots,h'_n)$ из $D(m,n)$ расстояние 
$d(c,c')=\sum^{m}_{i=1} d(s_i,s'_i)+\sum^{n}_{j=1}(s_j,s'_j)$.

\emph{Кодом} назовем произвольное подмножество вершин графа $D(m,n)$. Вершины, принадлежащие коду, будем называть \emph{кодовыми}. \emph{Кодовое расстояние} равно минимальному расстоянию между различными кодовыми вершинами. $(m+n,4^k,d)$ МДР кодом назовем код в графе $D(m,n)$ мощности $4^k$ для которого выполнено $d=2m+n-k+1$. 

Пусть дан набор координат $(i_1,\ldots,i_v;j_1,\ldots,j_w)$, где $1 \le i_1< \ldots <i_v \le m; 1 \le j_1 < \ldots < j_w \le n$. 

Для кода $C$ определим 
\emph{проекцию} $C_{i_1,\ldots,i_v;j_1,\ldots,j_w}$ и 
\emph{грань} $C^{a_1,\ldots,a_v;b_1,\ldots,b_w}_{i_1,\ldots,i_v;j_1,\ldots,j_w}$, где $a_i \in \VV Sh$ для $i=1,\ldots,v$ и $b_j \in \VV K$ для $j=1,\ldots,w$. 

Для начала определим:

$C_{i;}=\{(s_1,\ldots,s_{i-1},s_{i+1},\ldots,s_m;h_1,\ldots,h_n):
$ существует вершина $x \in \VV Sh  
$ такая, что $(s_1,\ldots,s_{i-1},x,s_{i+1},\ldots,s_m;h_1,\ldots,h_n) \in C\}$, $i \in \{1,\ldots,m\}$.

$C_{;i}=\{(s_1,\ldots,s_m;h_1,\ldots,h_{i-1},h_{i+1},\ldots,h_n): 
$ существует вершина $y \in \VV K 
$ такая, что $(s_1,\ldots,s_m;h_1,\ldots,h_{i-1},y,h_{i+1},\ldots,h_n) \in C\}$, $i \in \{1,\ldots,n\}$.
 
Тогда проекцию $C_{i_1,\ldots,i_v;j_1,\ldots,j_w}$ можно определить по правилам: 

$C_{i_1,\ldots,i_v;j_1,\ldots,j_w}$ $=(C_{i_1,\ldots,i_v;j_2,\ldots,j_{w}})_{;j_1}$ и $C_{i_1,\ldots,i_w;}$ $=(C_{i_2,\ldots,i_{w};})_{i_1;}$.

Также определим:

$C^{a;}_{i;}=\{(s_1,\ldots,s_{i-1},s_{i+1},\ldots,s_m;h_1,\ldots,h_n):(s_1,\ldots,s_{i-1},a,s_{i+1},\ldots,s_m;$ $h_1,\ldots,h_n) \in C\}$.

$C^{;b}_{;i}=\{(s_1,\ldots,s_m;h_1,\ldots,h_{i-1},h_{i+1},\ldots,h_n):(s_1,\ldots,s_m;h_1,\ldots,h_{i-1},b,$ $h_{i+1},\ldots,h_n) \in C\}$.

Тогда грань $C^{a_1,\ldots,a_v;b_1,\ldots,b_w}_{i_1,\ldots,i_v;j_1,\ldots,j_w}$ можно определить по следующим правилам:

$C^{a_1,\ldots,a_v;b_1,\ldots,b_w}_{i_1,\ldots,i_v;j_1,\ldots,j_w}$ $=(C^{a_1,\ldots,a_v;b_2,\ldots,b_{w}}_{i_1,\ldots,i_v;j_2,\ldots,j_{w}})^{;b_1}_{;j_1}$ и 
$C^{a_1,\ldots,a_v;}_{i_1,\ldots,i_v;}$ $=(C^{a_2,\ldots,a_{v};}_{i_2,\ldots,i_{v};})^{a_1;}_{i_1;}$.





Два кода $C$ и $C'$ назовем \emph{эквивалентными}, если существует 
набор автоморфизмов 
$\theta_1,\ldots,\theta_m$ графа $Sh$,
набор перестановок $\sigma_1,\ldots,\sigma_n$ из $S_4$
и перестановки координат $\tau_1  \in S_{m}$ и $\tau_2 \in S_{n}$ такие, что

$C'=\{(\theta_1 (s_{{\tau}_1(1)}),\ldots,\theta_m (s_{{\tau}_1(m)});
\sigma_1 (h_{{\tau}_2(1)}),\ldots,\sigma_n (h_{{\tau}_2 (n)})):
(s_1,\ldots,s_m;$ $ h_1,\ldots,h_n)
\in C\}$

Обозначим через $L_{m,n,k}$ число МДР кодов с параметрами $(m+n,4^k,d)$ с точностью до эквивалентности .

Граф $G=(V,E)$ называется сильно регулярным с параметрами $(v, k, $ $\lambda, \mu)$, если $G$ регулярный граф степени $k$ на $v$ вершинах, и любая пара смежных вершин имеет $\lambda$ общих соседей, и любая пара несмежных вершин имеет $\mu$ общих соседей. Графы $Sh$ и $K^2$ --- сильно регулярные с параметрами (16, 6, 2, 2).

В графе Шрикханде любая пара различных вершин $a$ и $b$ имеет ровно 2 общих соседа.
Обозначим их через $u(a,b)$ и $w(a,b)$.
Напомним, что вершины графа Шрикханде --- это элементы группы $Z^2_4$. 
Обозначим множества 
$A=\{01, 03, 10, 30, 11, 33\}$,
$B=\{02, 20, 22\}$, 
$C=\{12, 32, 13, 31, 21, 23\}$.
Элементы из $A$ и $C$ имеют порядок 4, элементы из $B$ имеют порядок 2.
 
\begin{lemman}\label{l:s}
Пусть $a$ и $b$ --- несмежные вершины графа Шрикханде.
Тогда если порядок элемента $(a-b)$ равен 2, то вершины $u(a,b)$ и $w(a,b)$ несмежны.
Если порядок элемента $(a-b)$ равен 4, то вершины $u(a,b)$ и $w(a,b)$ смежны.
\end{lemman}
\begin{proof}
Пусть элемент $(a-b)$ принадлежит множеству $B$.
Допустим $a=b+02$.
Тогда $u(a,b)=a+01$, $w(a,b)=a+03$, и эти вершины несмежны.
Если $a=b+20$ либо $a=b+22$, то доказывается аналогично.

Пусть теперь элемент $(a-b)$ принадлежит $C$.
Рассмотрим окрестность вершины $u(a,b)$. 
Тогда $a=u(a,b)+s_1$, $b=u(a,b)+s_2$, где $s_1, s_2 \in A$,
$s_1 \ne s_2$, вершины $s_1$ и $s_2$ несмежны, и $(s_1 -s_2) \in A$. 
Нетрудно убедиться, что тогда найдется $s \in A$, такой, что $s$ смежно с $s_1$ и $s_2$.
Тогда $w(a,b)=a+s$, и вершины $u(a,b)$ и $w(a,b)$ смежны.

\end{proof}

\begin{lemman}\label{l:p}
Для любой пары вершин $a, b$ графа Шрикханде и любого автоморфизма $\tau$ порядок элемента $(a-b)$
равен порядку элемента $(\tau(a)-\tau(b))$.
\end{lemman}
\begin{proof}
Если вершины $a$ и $b$ смежны, то вершины $\tau(a)$ и $\tau(b)$ также смежны, следовательно, $(a-b) \in A$ и $(\tau(a)-\tau(b)) \in A$.

Пусть $a$ и $b$ несмежны.
Предположим, что элементы $(a-b)$ и $(\tau(a)-\tau(b))$ имеют разный порядок.
Вершины $\tau(u(a,b))$, $\tau(w(a,b))$ будут общими соседями для вершин $\tau(a)$ и $\tau(b)$.
Тогда из леммы \ref{l:s} следует, что вершины $u(a,b)$ и $w(a,b)$ смежны тогда и только тогда, когда $\tau(u(a,b))$ и $\tau(w(a,b))$ несмежны. Противоречие с тем, что $\tau$ ---  автоморфизм.
\end{proof}

Максимальное независимое множество в графах $Sh$ и $K^2$ назовем \emph{кокликой}. Число независимоти графа Шрикханде равно 4.
Вершина $00$ из $Sh$ содержится в 4 различных кокликах:

$\{00,02,20,22\}$, $\{00,02,21,23\}$, $\{00,20,12,32\}$, $\{00,22,13,31\}$.

Третья и четвертая коклики эквивалентны второй (в качестве соответствующих автоморфизмов графа $Sh$ можно взять $\tau(ab)=ba$ и $\tau(ab)=(b-a)b$ соответственно).
Первая и вторая коклики неэквивалентны по лемме \ref{l:p}.
Коклики, содержащие вершину $s \in \VV Sh$, получаются прибавлением к каждому элементу из этих коклик элемента из $Z^2_4$, соответствующего вершине $s$.

С точностью до эквивалентности в графе $Sh$ ровно 2 различных коклики: $\{00,02,20,22\}$ и $\{00,02,21,23\}$. Первую назовем линейной, вторую полулинейной.
Эти коклики изображены на рис. \ref{f:mds}.
Все, с точностью до эквивалентности, разбиения графа Шрикханде на непересекающиеся коклики изображены на рис. \ref{f:latin}.

\begin{lemman}\label{l:koklika}
Пусть $G=(V,E)$ --- граф Шрикханде либо граф $K^2$, и 
$U=\{u_0, u_1, u_2, u_3\}$ --- произвольная коклика в $G$.
Тогда любая вершина из $V \backslash U$ смежна ровно с двумя вершинами из $U$.
\end{lemman}
\begin{proof}
Степень любой вершины графа $G$ равна 6. Следовательно количество ребер, таких, что один из концов ребра принадлежит $U$ равно 24. Мощность $|V \backslash U|$ равна 12. Предположим утверждение леммы неверно.
Тогда существует вершина $u \in V \backslash U$ смежная с 3 либо 4 вершинами из $U$. 

Рассмотрим окрестность вершины $u$. Докажем, что в этой окрестности не существует 3 или 4 попарно несмежных вершин, принадлежащих некоторой коклике. Это докажет лемму.

Если $G$ --- граф $K^2$, то окрестностью вершины $u$ будет объединение двух непересекающихся полных графов на 3 вершинах. Очевидно, что в этой окрестности нельзя выбрать 3 либо 4 попарно несмежных вершины.

Если $G$ --- граф Шрикханде, то окрестностью вершины $u$ будет цикл на 6 вершинах. Обозначим вершины цикла через  $a_1, a_2, a_3, a_4, a_5, a_6$, где вершины $a_i$ и $a_{i-1}$ смежны для $i=2,3,4,5,6$, а также смежны вершины $a_1$ и $a_6$. Очевидно, что в этой окрестности нельзя выбрать 4 попарно несмежных вершины. 
Предположим, что 3 вершины из окрестности принадлежат $U$. Тогда это либо вершины $u_1, u_3, u_5$, либо $u_2, u_4, u_6$. Без потери общности будем считать, что это вершины $a_1, a_3, a_5$. 
Обозначим через $W$ множество вершин, состоящее из вершины $u$ и вершин из ее окрестности.
Каждая вершина $a_i$, где $i=1,2,3,4,5,6$, смежна ровно с 3 вершинами из $W$, а следовательно и с 3 вершинами из $V \backslash W$. Легко увидеть, что 2 общих соседа вершин $a_1, a_3$ принадлежат множеству $W$ (вершины $u$ и $a_2$).
То же верно для вершин $a_3, a_5$, и для вершин $a_1, a_5$. Следовательно, есть ровно 9 вершин из $V \backslash W$ смежных с одной из вершин $a_1, a_3, a_5$.
Но мощность $|V \backslash W|=9$, 
следовательно, 
множество $\{a_1, a_3, a_5\}$ не может содержаться в какой-либо коклике. Лемма доказана.
\end{proof}

\begin{figure}
\centering
(a) \raisebox{-12mm}{\sh{nz nz nz nz}{zz nz zz nz}{nz nz nz nz}{zz nz zz nz}}\ 
(b) \raisebox{-12mm}{\sh{nz nz nz nz}{nz zz nz zz}{nz nz nz nz}{zz nz zz nz}}
\caption{Все коклики в $Sh$ с точностью до эквивалентности}
\label{f:mds}
\end{figure}

\begin{figure}
\centering
(a) \raisebox{-12mm}{\sh{nz vz nz vz}{zz wz zz wz}{nz vz nz vz}{zz wz zz wz}}\ 
(b) \raisebox{-12mm}{\sh{nz vz nz wz}{zz wz zz vz}{nz vz nz wz}{zz wz zz vz}}\ 
(c) \raisebox{-12mm}{\sh{nz vz zz wz}{zz wz nz vz}{nz vz zz wz}{zz wz nz vz}}
\caption{Все разбиения графа $Sh$ на непересекающиеся коклики с точностью до эквивалентности}
\label{f:latin}
\end{figure}


\begin{lemman}\label{l:kk}
Пусть $C$ --- $(2m+n, 4^k, 2m+n-k+1)$ МДР код. Если $2v+w=2m+n-k$, то в $C_{i_1,\ldots,i_v;j_1,\ldots,j_w}$ 
каждая вершина графа $D(m-v, n-w)$ встречается ровно один раз.
\end{lemman}
\begin{proof}
Если какая-то вершина в графе $D(m-v, n-w)$ встретится дважды, то расстояние между сооветствующими словами кода $C$
будет не больше $2m+n-k$, и мы получаем противоречие с кодовым расстоянием. 
С другой стороны число $4^k$ кодовых вершин  равно количеству вершин в графе $D(m-v, n-w)$.   
\end{proof}
Из леммы \ref{l:kk} следует, что если в $(m+n,4^k,d)$ МДР коде мы зафиксируем такой набор координат $i_1,\ldots,i_v;j_1,\ldots,j_w$, что $2v+w=k$, то значения в любой другой координате можно представить в виде функции, определенной на всем множестве вершин графа $D(v,w)$.

\begin{lemman}\label{l:prgr}
Пусть $C$ --- $(m+n, 4^k, d)$ МДР код
и $2v +w \le k$. Тогда:
\begin{enumerate}
\item Множество $C_{i_1,\ldots, i_v; j_1,\ldots, j_w}$ --- 
$((m-v)+(n-w), 4^{k}, d-2v-w)$ МДР код;

\item Множество ${C^{a_1,\ldots,a_{v}; b_1,\ldots, b_{w}}_{i_1,\ldots,i_{v}; j_1,\ldots, j_{w}}}$ 
--- $((m-v)+(n-w), 4^{k-2v-w}, d)$ МДР код.
\end{enumerate}
\end{lemman}


















\begin{lemman}\label{l:2ng}
Пусть даны два графа $G_1(V,E_1)$ и $G_2(V,E_2)$, 
каждый из которых явяется либо графом $Sh$, либо графом $K^2$, 
и $E_1 \cap E_2=\emptyset$. 
Тогда граф $G_3=(V, E_3=\overline{E_1 \cup E_2})$ --- объединение 4 непересекающихся графов $K_4$.
\end{lemman}
\begin{proof}
Графы $G_1$ и $G_2$ сильно регулярные с параметрами 
$(v,k,\lambda, \mu)=$ $(16,6,2,2)$ 
(т. е. регулярные графы нa $v$ вершинах степени $k$ такие, что любая пара смежных вершин имеет ровно $\lambda$ общих соседей, и любая пара несмежных вершин имеет ровно $\mu$ общих соседей).
Дополнение графа $\overline{G_1}$ сильно регулярный граф 
с параметрами $(16,9,4,6)$.

Определим для каждой вершины  $a$ из $V$ значения $S(a)$ и $N(a)$.
Обозначим через $V^2_{a}$ и $V^3_{a}$ множества вершин из окрестности вершины $a$ 
в графах $G_2$ и $G_3$ соответственно.
Тогда 

$$S(a)=|E_3 \cap (V^2_a \times V^2_a)|,$$
$$N(a)=|E_3 \cap (V^3_a \times V^3_a)|.$$

Так как $G_3$ --- регулярный граф степени 3, то $N(a)\le 3$.

Для произвольной вершины $a$ рассмотрим окресность этой вершины в графе $\overline{G_1}$ --- 
граф $(V^2_a \cup V^3_a,$ $(E_2 \cup E_3) \cap $ $((V^2_a \cup V^3_a)\times (V^2_a \cup V^3_a))$.
Это регулярный граф степени 4 (действительно, если взять любую вершину $b$ из окрестности $a$ в графе $\overline{G_1}$, 
то у нее с вершиной $a$ будет ровно $\lambda =4$ общих соседей).
Тогда в этом графе 18 ребер. 

С другой стороны число ребер равно ($S(a)$+6) (ребра, соединяющие вершины из $V^2_a$)
+ ($3\cdot 4$ $-N(a)$) (ребра, у которых хотя бы один из концов лежит в $V^3_a$).
Тогда $18=S(a)+6+12-N(a)$, и $S(a)=N(a)$.    

Посчитаем сумму $\sum_{a\in V} S(a)$.
Так как граф $G_2$ сильно регулярный с параметрами 
(16,6,2,2), 
то каждая пара вершин, соединенная  ребром из $E_3$,
содержится в окрестности ровно двух вершин в графе $G_2$, а следовательно 
каждое ребро из $E_3$ будет посчитано ровно два раза.
Следовательно, $\sum_{a \in V} S(a)=48$. 

А следовательно, $\sum_{a \in V} N(a)=48$,
а так как для любой вершины $a$ выполняется $N(a)\le 3$, то
$N(a)=3$ для любой вершины $a$, откуда и следует утверждение леммы.
\end{proof}

Обозначим через $S_m$ число таких четверок $(a, b, c, d)$, что $a, b, c, d$ --- целые неотрицательные числа, $a+b+c+d=m$ и $b \le c \le d$. 

\begin{lemman}\label{l:chislotroek}
$$S_m=m^3/36+7m^2/24+11m/12+1-(m \mod 2 )/8 - (m \mod 3 )/9.$$


\end{lemman}  

\begin{proof}
Будем рассматривать только такие четверки $(a,$ $b,c,d)$, что $a,b,c,d$ --- целые неотрицательные числа и $a+b+c+d=m$.

Представим $S_m$ в виде $S_m=M_m+N_m+K_m$, где 

$M_m$ --- число таких четверок $(a, b, c, d)$, что $b=c=d$;

$N_m$ --- число таких четверок $(a, b, c, d)$, что либо $b=c < d$, либо $b < c=d$;

$K_m$ --- число таких четверок $(a, b, c, d)$, что $b < c< d$.

Посчитаем $M_m$. По значению $b$ значение $a$ воостанавливается однозначно. Тогда так как число $b$ не превосходит $\lfloor m/3 \rfloor$, то  
$M_m=\lfloor m/3 \rfloor +1$.

Посчитаем $N_m$. Для начала посчитаем число $L_m$ таких четверок $(a, b,$ $ c, d)$, что 
$b \ne  c=d$. Для фиксированного числа $c$ число четверок $(a,b,c,c)$ равно $m-2c+1$. 
Число $c$ пробегает значения от 0 до $\lfloor m/2 \rfloor$. Тогда 
$L_m=(\sum^{i=\lfloor m/2 \rfloor}_{i=0}$ $(m-2i+1)$) $-M_m$ $=(((m+1)+(m-2 \lfloor m/2 \rfloor+1))/2) (\lfloor (m+2)/2) \rfloor)-M_m$=$\lfloor (m+3)/2 \rfloor \lfloor (m+2)/2 \rfloor$.
Здесь мы использовали равенство $(m+1-\lfloor m/2 \rfloor)=\lfloor (m+3)/2 \rfloor$. В его верности легко убедиться, рассмотрев случаи $m=2k$ и $m=2k+1$.
Число четверок $(a,b,c,d)$ таких, что $b=c \ne d$
также равно $L_m$.
Так как при подсчете $N_m$ из двух четверок $(a,b,b,c)$ и $(a,c,b,b)$, где $b \ne c$, мы выбираем ровно одну, то 
тогда $N_m=(L_m+L_m)/2=L_m$.
Тогда $$N_m=(\lfloor m/2 \rfloor+1) (\lfloor (m+1)/2 \rfloor+1)-M_m$$

Посчитаем $K_m$. Для начала посчитаем число таких четверок $(a,b,c,d)$, что числа $b, c, d$ попарно различны. Это число равно $(m+3)(m+2)(m+1)/6$
$-3L_m$ (число четверок, таких что среди $b, c, d$ ровно два различных числа) $-M_m$ (число четверок таких, что числа $b, c, d$ совпадают). Так как по условию $b < c < d$, то чтобы получить $K_m$ это число нужно поделить на 6 (чтобы упорядочить числа $b, c, d$). 
Тогда 
$$K_m=((m+3)(m+2)(m+1)/6-3L_m-M_m)/6.$$

Подставив выражения для $M_m, N_m, K_m$, получим, что $$S_m=(m+3)(m+2)(m+1)/36 
+(\lfloor m/2 \rfloor+1) (\lfloor (m+1)/2 \rfloor+1)/2
+(\lfloor m/3 \rfloor +1)/3.$$

Используем следующие равенства:
$$\lfloor m/2 \rfloor + \lfloor (m+1)/2 \rfloor =m,$$
$$\lfloor m/2 \rfloor \lfloor (m+1)/2 \rfloor = m^2/4-(m \mod 2)/4,$$
$$\lfloor m/3 \rfloor =m/3 -(m \mod 3)/3.$$
В верности первых двух равенств легко убедиться, рассмотрев случаи $m=2k$, $m=2k+1$, а в верности третьего --- рассмотрев случаи $m=3k$, $m=3k+1$, $m=3k+2$. 

Тогда подставив эти равенства в выражение для $S_m$, получим утверждение леммы.

\end{proof}

\section{Основная теорема.}\label{s:main}

Сформулиуем основную теорему.

\begin{theoreman}\label{t:main}
\begin{enumerate}
\item $L_{m,n,1}=m^3/36+7m^2/24+11m/12+1-(m \mod 2)/8 - ( m \mod 3 ) /9$.
\item При $4 \le 2m+n \le 6$ и $3 \le d \le 4$ значения для $L_{m,n,2m+n-d+1}$ представлены в таблице 1.
\item Если $2m+n=6$, то $L_{m,n,2}=0$.
\item Если $2m+n>6$ и $2<d<2m+n$, то $L_{m,n,k}=0$.
\end{enumerate}
\end{theoreman}

\begin{table}[h]
\begin{center}
\begin{tabular}{|c|c|c|c|c|c|c|c|}
\hline
$(m,n)$ & (2,0) & (1,2) & (2,1) & (1,3) & (2,2) & (1,4) & (3,0) \\
\hline
$d=3$   &   2   &   1   &   2   &   1   &   0   &   0   &    0 \\
\hline 
$d=4$   &   4   &   2   &   2   &   1   &   1   &   0   &    0 \\
\hline
\end{tabular}
\end{center}
\caption{Таблица значений $L_{m,n,2m+n-d+1}$}
\label{tab:latkv}
\end{table} 

Доказательство разобьем на несколько частей.


\section{$(m+n,4^1,2m+n)$ --- МДР коды}\label{s:mn1}

\begin{predln}\label{p:mn1}
$$L_{m,n,1}=m^3/36+7m^2/24+11m/12+1-(m \mod 2 )/8 - (m \mod 3)/9.$$
\end{predln}


Для начала докажем две следующие леммы:

\begin{lemman}\label{l:dhkvd}
Для любого $n \ge 0$ выполняется $L_{m,n,1}=L_{m,0,1}$.
\end{lemman}
\begin{proof}
Пусть даны два $(m+(i+1), 4^1, 2m+i+1)$ МДР кода $C$ и $C'$.
Легко видеть, что $C$ и $C'$ эквивалентны тогда и только тогда, когда проекции
$C_{;1}$ и $C'_{;1}$ эквивалентны.
Таким образом $L_{m,i,1}=L_{m,i+1,1}$,
откуда следует утверждение леммы.   
\end{proof}

\begin{lemman}\label{l:nm}
Пусть $U=\{u_0,u_1,u_2,u_3\}$ --- коклика в $Sh$.
Тогда:
\begin{enumerate}
\item если $U$ --- линейная, то существует автоморфизм $\tau$ графа $Sh$ такой, что:

$\tau(u_0)=00$, $\tau(u_1)=02$, $\tau(u_2)=20$, $\tau(u_3)=22$;
\item если $U$ --- полулиненая, то выполняется ровно одно из следующих утверждений:

1) существует автоморфизм $\tau$ графа $Sh$ такой, что 

$\tau(u_0)=00, \tau(u_1)=02, \tau(u_2)=21, \tau(u_3)=23$;

2) существует автоморфизм $\tau$ графа $Sh$ такой, что 

$\tau(u_0)=00$, $\tau(u_1)=21$, $\tau(u_2)=02$, $\tau(u_3)=23$; 

3) существует автоморфизм $\tau$ графа $Sh$ такой, что 

$\tau(u_0)=00$, $\tau(u_1)=21$, $\tau(u_2)=23$, $\tau(u_3)=02$.
\end{enumerate}
\end{lemman}
\begin{proof}
\begin{enumerate}

\item Пусть $U$ --- линейная. Тогда вершину $u_0$ можно перевести в вершину $00$ автоморфизмом $\theta(s)=s-u_0$. Тогда $\{u_1, u_2, u_3\}$ перейдет во множество $\{00,02,20,22\}$, а
все различные перестановки элементов множества $\{02, 20, 22\}$ дают, например, автоморфизмы $\tau_1(ab)=ab$, $\tau_2(ab)=ba$, $\tau_3(ab)=a(a-b)$, $\tau_4(ab)=(b-a)b$, $\tau_5(ab)=(a-b)a$, $\tau_6(ab)=(a-b)a$.

\item Пусть $U$ --- полулинейная. Для начала докажем то, что хотя бы одно из утверждений верно.  
$U$ некоторым автоморфизмов можно перевести во множество $\{00,02,21,23\}$. 
Тогда вершину $u_0$ можно перевести в вершину $00$ автоморфизмом $\theta(s)=s-u_0$. Тогда, так как $u_0$ принадлежит $\{02,21,23\}$ элемент $(-s)$ также принадлежит этому множеству. Тогда легко убедиться, что $\theta$ переводит $\{02,21,23\}$ в себя.
Поменять местами вершины $21$ и $23$ можно, например, автоморфизмом $\sigma(ab)=(-a)(-b)$. 
То, что верно ровно одно из утверждений следует из леммы \ref{l:p}. 

\end{enumerate}
\end{proof}

\begin{proof}[Предложения \ref{p:mn1}]
Из леммы \ref{l:dhkvd} следует, что достаточно рассматривать только МДР коды, в параметрах которых $n=0$.

Пусть $C$ --- $(m+0, 4^1, 2m)$ МДР код. 
После некоторого упорядочивания обозначим его вершины через $c^i=(s^i_1,\ldots,s^i_m)$, где $i=0, 1, 2, 3$ .
Также обозначим через $A_i=\{s^0_i,s^1_i,s^2_i,s^3_i\}$, $i=1,\ldots,m$.
Из кодового расстояния следует, что это будет коклика.

Назовем $(m+0, 4^1, 2m)$ МДР код $C$ \emph{приведенным}, если для некоторых $l, j, t$ таких, что $l \le j \le t$,  выполнено:

$s^0_i=00, s^1_i=02, s^2_i=21, s^3_i=23$, $i=1,\ldots,l$;

$s^0_i=00, s^1_i=21, s^2_i=02, s^3_i=23$, $i=l+1,\ldots,l+j$;

$s^0_i=00, s^1_i=21, s^2_i=23, s^3_i=02$, $i=l+j+1,\ldots,l+j+t$;

$s^0_i=00, s^1_i=02, s^2_i=20, s^3_i=22$, $i=l+j+t+1,\ldots,m$.

Так как для любой такой тройки $(l, j, t)$ приведенный код определяется единственным образом, обозначим этот код через $C^{l, j, t}$.

Обозначим через $T_m$ количество различных приведенных $(m+0,4^1,2m)$ МДР кодов. 
Это число равно количеству таких четверок $(m-l-j-t,l,j,t)$, что $l \le j \le t$.
Тогда из леммы \ref{l:chislotroek} следует, что 
$T_m=S_m$  .

Для $(m+0,4^1,2m)$ МДР кода $C$ определим четверку чисел $(a, b, c, d)_C$, где $b \le c \le d $.

Значение $a$ --- число линейных коклик во множестве $\{A_i:i=1,\ldots,m\}$.
Пусть $N_i(C)$ --- число элементов порядка 2 во множестве $\{(s^0_j-s^i_j):j=1,\ldots,m\}$, $i=1,2,3$.
Пусть $M_i(C)=N_i(C)-a$.
Упорядочив значения $M_1(C), M_2(C), M_3(C)$ 
мы получим значения $b, c, d$.

Докажем, что четверка $(a,b,c,d)_C$ не зависит от упорядочивания вершин кода.
Для значения $a$ это тривиально.
Достаточно доказать, что 
если поменять местами кодовые вершины $c^1, c^2, c^3$, либо поменять местами вершины $c^0$ и $c^1$, то значения $b, c, d$ не изменятся.
Если поменять местами кодовые вершины $c^1, c^2, c^3$, то изменится только порядок $N_1(C), N_2(C), N_3(C)$, следовательно, значения $b, c, d$ не изменятся.
Определим
$N'_i(C)$ --- количество элементов порядка 2 во множестве $\{(s^1_j-s^i_j:j=1,\ldots,m)\}$, $i=0,2,3$.
Очевидно, что $N_1(C)=N'_0(C)$.
Также $N_2(C)=N'_3(C)$. Действительно, так как для любого $j=1,\ldots,m$ множество $\{s^0_j, s^1_j, s^2_j, s^3_j\}$ будет кокликой, то порядок элемента $(s^1_j-s^3_j)$ равен порядку элемента $(s^0_j-s^2_j)$
(достаточно проверить это для коклик $\{00,02,20,22\}$ и $\{00,02,21,23\}$).
Аналогично $N_3(C)=N'_2(C)$.
Следовательно, если поменять местами вершины $c^0$ и $c^1$, то
значения $b, c, d$ не изменятся.

Если два $(m+0,4^1,2m)$ МДР кода $C$ и $C'$ эквивалентны, то четверки $(a,b,c,d)_C$ и $(a,b,c,d)_{C'}$ совпадают.
Действительно, очевидно, что перестановка координат в коде $C$ не меняет $(a,b,c,d)_C$. Также по лемме \ref{l:p} набор автоморфизмов 
графа Шрикханде не меняет $(a,b,c,d)_C$.

Для приведенного кода $D=C^{l,j,t}$ четверка $(a,b,c,d)_D$ равна $(m-l-j-t, l, j, t)$. 
Тогда, если тройки чисел $(l_1, j_1, t_1)$ и $(l_2, j_2, t_2)$ не совпадают, то коды $C^{l_1, j_1, t_1}$ и $C^{l_2, j_2, t_2}$ не эквивалентны.
Тогда все приведенные коды попарно неэквивалентны,
следовательно, $T_m \le L_{m,n,1}$.

Докажем, что произвольный $(m+0,4^1,2m)$ МДР код $C$ эквивалентен некоторому приведенному коду.
Этому коду соответствует четверка $(a,b,c,d)_C$.
Тогда можно взять перестановку координат $\pi$ такую, что для кода $C'$, полученного из кода $C$ перестановкой координат $\pi$, выполнено:   

элемент $(s^0_j-s^1_j)$ имеет порядок 2, $j=1,\ldots,b$;

элемент $(s^0_j-s^2_j)$ имеет порядок 2, $j=b+1,\ldots,b+c$;

элемент $(s^0_j-s^3_j)$ имеет порядок 2, $j=b+c+1,\ldots,b+c+d$;

множество $\{s^0_j, s^1_j, s^2_j, s^3_j\}$ --- линейная коклика, $j=b+c+d+1,\ldots,m$. 

Тогда по лемме \ref{l:nm} существует набор автоморфизмов $\tau_1,\ldots,\tau_m$ графа Шрикханде, переводящий код $C'$ в приведенный код $C^{b,c,d}$.

Отсюда следует, что $L_{m,n,1} \le T_m$.
Тогда $L_{m,n,1}=T_m$, что доказывает предложение.

\end{proof}


\section{$(m+n, 4^2, 2m+n-1)$ --- МДР коды}\label{s:mn2}


\begin{predln}\label{p:mn2}
\begin{enumerate}
\item $L_{2,0,2}=2$, $L_{1,2,2}=1$.
\item $L_{2,1,2}=2$, $L_{1,3,2}=1$.
\end{enumerate}
\end{predln}
\begin{proof}
Для начала рассмотрим МДР коды с параметрами $(2+0, 4^2, 3)$ и $(1+2, 4^2, 3)$.
Пусть $C$ --- один из таких кодов.
Представим кодовые слова в виде $(a,f(a))$, где $a \in \VV Sh$,
а $f(a) \in \VV Sh$
в случае $(2+0,4^2,3)$ МДР кода, 
$f(a) \in \VV K^2$ в случае $(1+2,4^2,3)$ МДР кода. Из кодового расстояния следует, что отображение $f$ будет биекцией.

Для кода $C$ определим разбиение $L^C$ графа Шрикханде на коклики. 

Для начала определим  графы $G_1(V_1,E_1)$, $G_2(V_2, E_2)$ и $G_f(V_1,E_f)$, где $V_1=\VV Sh, V_2=\{f(a):a \in V_1\}$ и 

$$E_1=\{(a_1,a_2):d(a_1,a_2)=1, a_1, a_2 \in V_1\},$$ 
$$E_2=\{(f(a_1),f(a_2)):d(f(a_1),f(a_2))=1, f(a_1), f(a_2) \in V_2\},$$
$$E_f=\{(a_1,a_2):d(f(a_1),f(a_2))=1, a_1, a_2 \in V_1\}.$$

Очевидно, $G_1$ --- граф Шрикханде, граф $G_2$ ---граф Шрикханде либо граф $K^2$.
Граф $G_f$ изоморфен графу $G_2$, следовательно, также является графом Шрикханде либо графом $K^2$.

Множества ребер $E$ и $E_f$ не пересекаются (если некоторое ребро $(a,b)$ лежит в пересечении,
то между кодовыми вершинами $(a, f(a))$ и $(b, f(b))$ расстояние равно 2).


Из леммы \ref{l:2ng} следует, что существует единственное разбиение $V=L_0\cup L_1 \cup L_2 \cup L_3$, такое, что : 
$|L_i|=4$ для каждого $i=0, 1 , 2, 3$, 
и вершины $a$ и $b$ графа $G_3=(V_1, \overline{E_1 \cup E_f})$ смежны тогда и только тогда, когда $a$ и $b$ принадлежат одному и тому же множеству $L_i$ для некоторого $i=0,1,2,3$.
Множества вершин $L_0, L_1, L_2 ,L_3$ будут попарно непересекающимися кокликами в графах $G$ и $G_f$.
Тогда определим $L^C=\{L_0, L_1, L_2, L_3\}$.

Для каждой вершины $a$ графа Шрикханде обозначим через $L^{C}(a)$ ту коклику разбиения, которой она принадлежит. 

Из построения следует:

\begin{lemman}\label{l:d20d12}
Пусть $C$ --- МДР код с параматрами $(2+0, 4^2, 3)$ либо $(1+2, 4^2, 3)$. 
Тогда расстояние между двумя различными кодовыми вершинами 
$(a, f(a))$ и $(b, f(b))$ равно 4
тогда и только тогда, когда $b \in L^C(a)$, и равно 3 иначе. 
\end{lemman}






Также определим $R_i=\{f(a):a \in L_i\}$. 
Для каждой вершины $a$ обозначим через $R^{C}(a)$ то множество из этого разбиения, которому оно принадлежит.
По определению вершина $a \in L^{C}(b)$ тогда и только тогда, когда вершина $f(a) \in R^{C}(f(b))$.
Также по лемме \ref{l:d20d12} $R_i$ будет кокликой в графе $G_2$ для каждого $i=0,1,2,3$.

Пусть $L$ --- разбиение графа Шрикханде на коклики. Тогда если вершины $a$ и $b$ принадлежат одной коклике, то обозначим это через $aLb$, и обозначим $a\overline{L}b$ иначе. 

Для разбиения $L$ единственным образом можно определить граф $G_L=(V, E_L)$, где $V=\VV Sh$ и
$$E_L=\{(a,b): d(a,b)=2, a \overline{L} b\}.$$

Пусть $L^C=L$ для некоторого кода $C$. Тогда $f$ такое отображение, что $d(f(a),f(b))=1$ тогда и только тогда, когда $(a,b) \in E_L$.
Отсюда следует, что если $L^{C}=L^{C'}$ для некоторого кода $C'$, то $d(f(a),f(b))=d(f'(a),f'(b))$ для любой пары вершин из $Sh$, следовательно, коды $C$ и $C'$ эквивалентны. 
Также если разбиения $L^{C}$ и $L^{C'}$ эквивалентны, то и коды $C$ и $C'$ эквивалентны.
Действительно, так как разбиения эквалентны, то код $C'$ эквивалентен некоторому коду $C''$, такому, что $L^{C''}=L^{C}$. Тогда коды $C$ и $C''$ эквивалентны, а следовательно, коды $C$ и $C'$ также эквивалентны.
Тогда если для некоторого разбиения $L$ найти такой код $C$, что для отображения $f$ из $V_1$ в $V_2$ будет $E_L=\{(a,b):d(f(a),f(b))\}$, то любой другой код $C'$ такой, что разбиение $L^{C'}$ эквивалентно разбиению $L^{C}$, эквивалентен $C$.

В графе Шрикханде существует три разбиения на коклики с точностью до эквивалентности. Эти разбиение изображены на рис. \ref{f:latin}.

Если $L$ --- разбиение а), то граф $G_L$ является графом Шрикханде. 
Код, соответствующий этому разбиению приведен в Приложении в таблице \ref{tab:2021}.
При этом разбиение $\{R_i:i=0,1,2,3\}$, где $R_i=\{f(a):a \in L_i\}$, будет также разбиением а).

Если $L$ --- разбиение б), то граф $G_L$ является графом Шрикханде.
Код, соотвествующий этому разбиению приведен в Приложении в таблице \ref{tab:2022}.
При этом разбиение $\{R_i:i=0,1,2,3\}$, где $R_i=\{f(a):L_i\}$, будет разбиением б).
Тогда МДР коды, соответствующие разбиению а) и б) неэквивлентны.

Если $L$ --- разбиение в), то граф $G_L$ является графом $K^2$.
Код, соответствующий этому разбиению приведен в Приложении в таблице \ref{tab:1221}.

Таким образом мы получаем, что
$(2+0, 4^2, 3)$ МДР кодов будет 2 с точностью до эквивалентности, 
а $(1+2,4^2,3)$ 
МДР код будет единственен с точностью до эквивалентности.
Это доказывает пункт 1 предложения.

 

Пусть теперь $C$ --- МДР код с параметрами $(2+1, 4^2, 4)$ ( с параметрами $(1+3, 4^2, 4)$).
Обоначим его вершины через $(a, f(a), g(a))$, где $a \in Sh$, $f(a) \in \VV Sh$ ( $f(a) \in \VV K^2$) , $g(a)\in \VV K$.
Проекция $D=C_{;1}$ --- МДР код с параметрами $(2+0, 4^2, 3)$ (для $(1+3,4^2,4)$ МДР кода $D=C_{;3}$ ---  $(1+2, 4^2, 3)$ МДР код).
Из леммы \ref{l:d20d12} следует, что если $b$ не принадлежит $L^{D}(a)$, то из кодового расстояния получаем, что $g(a) \ne g(b)$. Так как для каждого $i=0,1,2,3$ выполнено
$|C^{;i}_{;1}|=4$ ($C^{;i}_{;3}$), то $g(a)=g(b)$ тогда и только тогда, когда $a \in L^{D}(b)$, и по заданному коду $C_{;1}$ код $C$ восстанавливается однозначно с точностью до перестановки из $S_4$. 
Отсюда следует, что $L_{2,1,2}=L_{2,0,2}$ и $L_{1,3,2}=L_{1,2,2}$.
\end{proof} 

Также из построения следует:


\begin{lemman}\label{l:d21d13}
Пусть $C$ --- МДР код с параметрами $(2+1, 4^2, 4)$ либо $(1+3, 4^2, 4)$.
Тогда расстояние между любыми двумя различными кодовыми вершинами равно 4.
\end{lemman}



\section{$(m+n, 4^3, 2m+n-2)$ МДР коды}\label{s:mn33}

\subsection{МДР коды с параметрами $(2+1, 4^3, 3)$ и $(1+3, 4^3, 3)$}\label{s:213}
\begin{predln}\label{p:mn3}
$L_{2,1,3}=2$, $L_{1,3,3}=1.$
\end{predln}

В данном подразделе будем рассматривать МДР коды с параметрами $(2+1, 4^3, 3)$ и $(1+3, 4^3, 3)$.
Пусть $C$ --- один из таких кодов.
Представим кодовые вершины  
в виде 
$(f_k(a),a,k)$, где $k \in \VV K$, $f_k(a) \in \VV Sh$ и  
$a$ --- в первом случае вершина графа Шрикханде, во втором вершина графа $K^2$.

Обозначим через $D_i$ грань $C^{;i}_{;1}$ ($D_i=C^{;i}_{;3}$ в случае $(1+3,4^3,3)$ МДР кода). 

Для начала докажем следующую лемму: 

\begin{lemman}\label{l:qwer}
Пусть $C$ --- МДР код с параметрами $(2+1, 4^3, 3)$ либо $(1+3,4^3,3)$. Тогда:
\begin{enumerate}
\item[\rm(i)] для любой вершины $a$ 
и любых различных $i,j\in\{0,1,2,3\}$
выполняется $d(f_i(a),f_j(a))=2$;
\item[\rm(ii)] для произвольной вершины $a$ и для любого $i=1,2,3$
$$\{f_k(a):k=0,1,2,3\}=
L^{D_0}(f_0(a))=
L^{D_i}(f_i(a));$$
\item[\rm(iii)] для любого $i=0,1,2,3$ 
$$R^{D_i}(a) =R^{D_0}(a);$$
\item[\rm(iv)] для любого $i=0,1,2,3$ и любой пары вершин $a$ и $b$ 
$$d(f_0(a),f_0(b))=d(f_i(a),f_i(b)).$$
\end{enumerate}
\end{lemman}
\begin{proof}
\begin{enumerate}
\item[(i)]
Так как расстояние между соответствующими вершинами 
$(f_i(a),a,i)$ и $(f_j(a),a,j)$ кода $C$
не меньше трех, 
то расстояние между вершинами
$f_i(a)$ и $f_j(a)$ равно $2$.
\item[(ii)] 

По лемме~\ref{l:kk} найдется $b$
такое, что $(f_1(a),b,0) \in C$ 
(т.е. $f_0(b)=f_1(a)$).

Так как расстояние между вершинами
$(f_1(a),a,1)$ и $(f_1(a),b,0)$ 
из $C$ не меньше $3$,
расстояние между $a$ и $b$ равно $2$.
Тогда из (i) следует, что
расстояние между вершинами 
$(f_0(a),a,0)$ и $(f_1(a),b,0)$ равно $4$.
Следовательно, по лемме 
\ref{l:d20d12} имеем $b \in R^{D_0}(a)$ и, следовательно, 
$f_1(a)\in L^{D_0}(f_0(a))$.
Аналогично получаем
$f_2(a)\in L^{D_0}(f_0(a))$ и 
$f_3(a)\in L^{D_0}(f_0(a))$,
кроме того, тривиально 
$f_0(a)\in L^{D_0}(f_0(a))$.
В итоге имеем
$\{f_k(a):k=0,1,2,3\}=
L^{D_0}(f_0(a))$, оставшиеся три равенства аналогичны.


\item[(iii)]
Напомним, что 
$\{ L^{D_i}(b) : b\in \VV \mathrm{Sh} \}$ --- ассоциированное с кодом 
$D_i$ разбиение вершин графа Шрикханде на $4$ коклики. Как следует из (ii),
это разбиение не зависит от $i$.
Из симметрии между первой и второй координатой следует, что разбиение 
$\{ R^{D_i}(b) : b\in \VV \mathrm{Sh} \}$ также не зависит от $i$, откуда следует требуемое утверждение.


\item[(iv)] Пусть $a\ne b$. 
По лемме \ref{l:d20d12} расстояние $d((f_i(a),a,i),(f_i(b),b,i))$ 
равно $4$, 
если $b \in R^{D_i}(a)$,
и равно $3$ в противном случае.
Но согласно (iii) имеем 
$R^{D_i}(a)=R^{D_0}(a)$,
откуда следует, 
что это расстояние не зависит от $i$. 
А значит, и расстояние между 
$f_i(a)$ и $f_i(b)$ не зависит от $i$.
\end{enumerate}
\end{proof}

\begin{lemman}\label{l:autnm}
Пусть $U=\{U_0, U_1, U_2, U_3\}$ --- разбиение графа Шрикханде на непересекающиеся коклики. 
Тогда существует единственный набор из трех автоморфизмов графа Шрикханде ---
$\tau_1,\tau_2,\tau_3$ такой, что
для любого $j=0,1,2,3$, для любого $i=1,2,3$ и для любой вершины графа Шрикханде $s$:

1) если $s \in U_j$, то $\tau_i(s) \in U_j$;

2) $d(\tau_i(s),s)=2$;

3) $d(\tau_i(s),\tau_j(s))=2$ при $i \ne j$.

%
%
%
%
%
\end{lemman}

\begin{proof}


Обозначим $A=\{02,20,22\}$, $B=\{12,32,21,23,13,31\}$.

Для начала рассмотрим случай, когда $U$ --- разбиение a) из рис. \ref{f:latin}.
Пусть $\tau$ --- автоморфизм, удовлетворяющий условиям 1 и 2.
Тогда для любой вершины $s \in \VV Sh$ либо $\tau(s)=s+02$, либо $\tau(s)=s+20$, либо $\tau(s)=s+22$.

Пусть $\tau(a)=a+02$ для некоторой вершины $a$.
Докажем, что тогда $\tau(s)=s+02$ для любой вершины $s \in \VV Sh$.
Вершины $a+01$ и $a+03$ принадлежат одной и той же линейной коклике $U_i$ для некоторого $i=0,1,2,3$.
Так как $a$ переходит в вершину $a+02$, то вершины $a+01$ и $a+03$ переходят в вершины из $U_i$, принадлежащие окрестности вершины $a+02$. По лемме \ref{l:koklika} таких вершин ровно 2, и это вершины $a+01$ и $a+03$. Тогда по условию 2 $\tau(a+01)=a+03$ и $\tau(a+03)=a+01$. Так как у вершин $a+01$ и $a+03$ 2 общих соседа, то множество вершин $\{a,a+02\}$ переходит во множество $\{a,a+02\}$, а следовательно,  $\tau(a+02)=a$.
Предположим, что $\tau(b) \ne b+02$ для некоторой вершины $b \in \VV Sh$. Пусть $\tau(b)=b+20$. Тогда, аналогично, $\tau(b+10)=b+30$, $\tau(b+20)=b$, $\tau(b+30)=b+10$. 
Множества вершин $\{a, a+01, a+02, a+03\}$ и $\{b,b+10,b+20, b+30\}$ пересекаются в некоторой вершине $c$.
Тогда $\tau(c)=c+02$ и $\tau(c)=c+20$, и мы получаем противоречие. Аналогично, если $\tau(b)=b+22$.

Аналогично доказывается для остальных $x \in A$, что если $\tau(a)=a+x$, для некоторой вершины $a$, то $\tau(s)=s+x$ для любой вершины $s \in \VV Sh$.
Тогда набор автоморизмов $\tau_1(s)=s+02$, $\tau_2(s)=s+20$, $\tau_3(s)=s+22$ будет единственным, удовлетворяющим условиям 1-3.

Рассмотрим случай, когда $U$ --- разбиение б) или в) из рис. \ref{f:latin}.

Пусть $a, b, c$ --- различные попарно смежные вершины в $Sh$. Тогда если вершина $d$ не совпадает с  $a, b, c$ и смежна с двумя из этих вершин, скажем с $a$ и $b$, то по значениям $\tau(a), \tau(b), \tau(c)$ однозначно восстанавливается значение $\tau(d)$. Это следует из того, что $\tau(c)$ и $\tau(d)$ будут общими соседями вершин $\tau(a)$ и $\tau(b)$, а так как у любой пары вершин в $Sh$ ровно 2 общих соседа, то значение $\tau(d)$ определяется однозначно.
Тогда, используя это свойство, нетрудно убедиться, что по значениям $\tau(a), \tau(b), \tau(c)$ однозначно задаются значения $\tau(s)$ для любой вершины $s \in \VV Sh$.

Пусть $\tau$ --- автоморфизм, удовлетворяющий условиям 1 и 2.
Докажем, что если $\tau(s)=s+x$, где $x \in B$,
для некоторой вершины $s$, то $\tau$ определяется однозначно по разбиению $U$.
Вершины $s$ и $s+x$ имеют 2 общих соседа, скажем вершины $a$ и $b$.
Так как $x$ имеет порядок 4, то по лемме \ref{l:s} эти вершины смежны, а следовательно,
принадлежат разным кокликам из $U$, скажем $U_i$ и $U_j$, где $i \ne j$.
По лемме \ref{l:koklika} в окрестности вершины $s$ есть ровно одна вершина 
$c \ne a$, принадлежащая $U_i$, а также ровно одна вершина $d \ne b$, принадлежащая $U_j$. Тогда так как по лемме \ref{l:koklika} в окрестности $s+x$ ровно 2 вершины из $U_j$, вершина $a$ и некоторая вершина $e$, то тогда либо $\tau(a)=a$, либо $\tau(c)=a$. Но первое неверно по условию 2, следовательно, $\tau(c)=a$.
Аналогично, $\tau(d)=b$.
Тогда, так как вершины $a, b, s$ попарно смежны, то по значениям $\tau(a), \tau(b), \tau(s)$ автоморфизм $\tau$ восстанавливается однозначно.
Таким образом, если существует автоморфизм $\tau$, удовлетворяющий условиям 1 и 2, такой, что $\tau(s)=s+x$, где $x \in B$, то такой автоморфизм единственный.
В разбиении $U$ найдется полулинейная коклика. Тогда эту коклику можно представить в виде $\{s, s+x, s+y, s+z\}$, где $x, y \in A$, $z \in B$, $s \in \VV Sh$. 
Тогда для набора автоморфизмов $\tau_1, \tau_2, \tau_3$, удовлетворяющих условиям 1-3,
$\tau_1(s)=s+x$, $\tau_2(s)=s+y$, $\tau_3(s)=s+z$. Автоморфизмы $\tau_1$ и $\tau_2$ определяются однозначно. Тогда, так как для любой вершины $s \in U_i$ и любого $i=0,1,2,3$ множество $\{s, \tau_1(s), \tau_2(s), \tau_3(s)\}$ равно $U_i$, то для любой вершины $s$ значение $\tau_3(s)$ определяется однозначно по значениям $\tau_1(s)$ и $\tau_2(s)$.
Таким образом если существует набор автоморфизмов $\tau_1, \tau_2, \tau_3$, то он единственненый.

Осталось привести такой набор автоморфизмов для каждого разбиения.

Пусть $U_0=\{00,02,21,23\}$, $U_1=\{01,03,22,20\}$, $U_2=\{10,12,30,32\}$, $U_3=\{11,13,31,33\}$.
Тогда 

$\tau_1(ab)=a(a-b)+21$, $\tau_2(ab)=ab+02$, $\tau_3(ab)=a(a-b)+23$.

Пусть $U_0=\{00,02,21,23\}$, $U_1=\{01, 03, 20, 22\}$, $U_2=\{10,12,31,33\}$, $U_3=\{11,13,30,32\}$.
Тогда 

$\tau_1(ab)=ab+02$, $\tau_2(ab)=ab+21$, $\tau_3(ab)=ab+23$.

Пусть $U'=\{U_0, U'_1, U'_2, U'_3\}$ и $U=\{U_0, U_1, U_2, U_3\}$ эквивалентные разбиения.
Тогда существует автоморфизм $\phi$, переводящий $U$  в $U'$, и $U'_j=\{\phi(x):x\in U_j\}$, $j=0,1,2,3,$

Тогда если набор автоморфизмов $\tau_1, \tau_2, \tau_3$ удовлетворяет условиям 1-3 для разбиения $U$, то набор автоморфизмов $\sigma_i(s)=\phi(\tau_i(\phi^{-1}(s)))$, $i=1,2,3$, удовлетворяет условиям 1-3 для разбиения $U'$.

Докажем, что набор автоморфизмов $\sigma_i$, где $i=1,2,3$, удовлетворяет условиям 1-3.
Пусть $s \in U'_j$. Тогда $\phi^{-1}(s) \in U_j$, $\tau_i(\phi^{-1}(s)) \in U_j$, $\phi(\tau_i(\phi^{-1}(s))) \in U'_j$, т.е. условие 1 выполняется.
Если для некоторого $i=1,2,3$ автоморфизм $\sigma_i$ не удовлетворяет условию 2, то из условия 1 следует, что $\sigma_i(s)=s$ для некоторой вершины $s$.
Но тогда $\tau_i(\phi^{-1}(s))=\phi^{-1}(s)$, и мы получаем противоречие.
Если не выполняется условие 3, то для некоторых $i \ne j$, и некоторой вершины $s$,
$\sigma_i(s)=\sigma_j(s)$, а следовательно, $\tau_i(\phi^{-1}(s))=\tau_j(\phi^{-1}(s))$, и мы получаем противоречие.

Утверждение леммы доказано.
\end{proof}

\begin{proof}[Предложения \ref{p:mn3}]
Пусть $C$ --- МДР код с параметрами $(2+1,4^3,3)$ (либо $(1+3,4^3,3))$.
Грань $D_0=C^{;0}_{;1}$ ($D_0=C^{0;}_{3;}$) является  МДР кодом с параметрами $(2+0,4^2,3)$ ($(1+2,4^2,3)$ соответственно).
Обозначим через $L_0, L_1, L_2, L_3$ --- разбиение на коклики, определенное следующим образом: вершины
$a$ и $b$ принадлежат одному и тому же множеству тогда, и только тогда, когда $a \in L^{D_0}(b)$.
Так как $f_0$ --- биекция, то для $i=1,2,3$ можно обозначить $\tau_i(a)=f_i (f^{-1}_0)$.
Из пункта (iv) леммы \ref{l:qwer} следует, что $\tau_i$ --- автоморфизм для любого $i=1,2,3$.
По пункту (i) леммы \ref{l:qwer} для любой вершины $s=f_0(a)$ графа Шрикханде  и для любого $i=1,2,3$ выполнено $d(\tau_i(s),s)=2$.
Также по пункту (i) леммы \ref{l:qwer} для любой вершины $s=f_0(a)$ и любых $i \ne j$ выполнено $d(\tau_i(s),\tau_j(d))=2$.
Из пункта (iii) леммы \ref{l:qwer} следует, что для любой вершины $s=f_0(a)$ и любого $i=1,2,3$, и любого $j=0,1,2,3$ выполнено: если 
$s \in L_j $, то $\tau_i(s) \in L_j$.
Тогда из леммы \ref{l:autnm} следует, что если для любой вершины $a \in \VV Sh$ известно значение $f_0(a)$, то для любой вершины $a \in \VV Sh$ и любого $i=1,2,3$ значение $f_i(a)$ восстанавливается с точностью до перестановки значений $i$.
Из этого следует, что если эквивалентны  
грани $C^{;0}_{;1}$ и $C'^{;0}_{;1}$ ($C^{;0}_{;3}$ и $C'^{;0}_{;3}$ для $(1+3,4^3,3)$ МДР кода) МДР кодов $C$ и $C'$ с параметрами $(2+1,4^3,3)$ ($(1+3,4^3,3)$), то коды $C$ и $C'$ эквивалентны.

Так как $(2+0,4^2,3)$ МДР кодов 2 с точностью до эквивалентности, то $L_{2,1,3}=L_{2,0,2}=2$.
Аналогчино, $L_{1,3,3}=L_{1,2,2}=1$.
\end{proof}

\subsection{МДР коды с параметрами $(2+2,4^3,4)$ или $(1+4,4^3,4)$}\label{s:mn34}


\begin{predln}\label{p:224}
$L_{2,2,3}=1$, $l_{1,4,3}=0$.
\end{predln}
\begin{proof}
При доказательстве будем рассматривать $(2+2,4^3,4)$ МДР коды. 
Для $(1+4,4^3,4)$ МДР кодов рассуждения аналогичны.
Пусть $C$ --- МДР код с параметрами $(2+2,4^3,4)$ (с параметрами $(1+4,4^3,4)$).
Обозначим его вершины через $(f_i(a),a,i,g_i(a))$, где $f_i(a) \in \VV Sh$,
$a \in \VV Sh$ $(a \in \VV K^2)$, и 
$i, g_i(a) \in \VV K$.

Обозначим через $P=C_{;2}$ (для $(1+4,4^3,4)$ МДР кода $P=C_{;4}$).
Обозначим через $D_i=P^{;i}_{;1}$ ($D_i=P^{;i}_{;3}$ сответственно).
По лемме \ref{l:qwer} разбиения $\{L^{D_i}(a):a \in \VV Sh\}$ 
и $\{R^{D_i}(a):a \in \VV Sh\}$  не зависят от $i$.
Обозначим эти разбиения через $L$ и $R$ соответственно, а
коклики из этих разбиений через $L_i$ и $R_i$, $i=0,1,2,3$.
Если вершины $a,b$ принадлежат одному и тому же множеству, то обозначим это через $a \stackrel{L} \sim b$ и $a \stackrel{R} \sim b$ соответственно.

Если $C$ --- $(2+2,4^3,4)$ МДР код, то $P$ --- $(2+1,4^3,3)$ МДР код. Таких кодов 2 с точностью до эквивалентности. Эти коды приведены в Приложении в таблицах \ref{tab:2131} и \ref{tab:2132}.
Если $P$ эквивалентен коду из таблицы \ref{tab:2131}, то разбиение $L$ будет разбиением а) из рис. \ref{f:latin}.
Если $P$ эквивалентен коду из таблицы \ref{tab:2132}, то разбиение $L$ будет эквивалентно разбиению б) из рис. \ref{f:latin}.
Если $C$ --- $(1+4,4^3,4)$ МДР код, то $P$ --- $(1+3,4^3,3)$ МДР код. Тогда код $P$ эквивалентен коду из таблицы \ref{tab:1331}, а 
тогда $L$ эквивалентно разбиению в) из рис. \ref{f:latin}.

Докажем, что $L$ --- разбиение a) из рис. \ref{f:latin}. 
Из этого будет следовать, что МДР код $P$ эквивалентен $(2+1,4^3,3)$ МДР коду из таблицы \ref{tab:2131} из Приложения, а $(1+4,4^3,4)$ МДР кодов не существует.

Для начала докажем ряд утверждений:



1. \emph{Для любого $i=0,1,2,3$ выполняется $L_t$ $=\{f_i(a):a \in R_t\}.$}

По пунктам (ii) и (iii) леммы \ref{l:qwer}.

2. \emph{Для любой вершины $a$ если $i \ne j$, то $d((f_i(a),(f_j(a))=2.$}

По пункту (i) леммы \ref{l:qwer}.

3. \emph{Для любых $i \ne j$ и для любой вершины $a$ выполняется $g_i(a) \ne g_j(a).$}

Иначе получим противоречие с кодовым расстоянием.

4. \emph{Для любого $i=0,1,2,3$ равенство  
$g_i(a)=g_i(b)$ имеет место тогда и только тогда, когда $a\stackrel{R}\sim b.$}

Рассмотрим грань $U=C^{;i}_{;1}$ $(C^{;i}_{;3})$.
По лемме \ref{l:d21d13} расстояние между вершинами 
$(f_i(a), a, g_i(a))$ и $(f_i(b), b, g_i(b))$ равно $4$.
Получается, что $g_i(a)=g_i(b)$ 
тогда и только тогда, когда
расстояние между $(f_i(a), a)$ и $(f_i(b), b)$ равно $4$.
По лемме \ref{l:d20d12} для проекции $U_{;1}$ $(U_{;3})$ последнее эквивалентно соотношению 
$a \stackrel{R}\sim b$. 

5. \emph{Если $i \ne 0$ и $g_0(a)=g_i(b)$, то $d(f_0(a),f_0(b))$ $=d(f_0(a),f_i(b)).$}

Из пункта 3 следует, что $g_0(b) \ne g_i(b)$, следовательно, $g_0(a) \ne g_0(b)$.
Из пункта 4 следует, что $a \not\stackrel{R}\sim b$.  
Вершины $(f_0(a), a, g_0(a))$ и $(f_0(b), b, g_0(b))$ 
принадлежат грани $C^{;0}_{;1}$ $(C^{;0}_{;3})$, следовательно, 
по лемме \ref{l:d21d13} расстояние между ними равно 4, и $d((f_0(a), a), (f_0(b), b))=3$.
Вершины $(f_0(a), a, 0)$ и $(f_i(b), b, i)$ 
принадлежат грани $C^{;g_0(a)}_{;2}$ $(C^{;g_0(a)}_{;4})$, следовательно, 
по лемме \ref{l:d21d13} расстояние между ними равно 4, и $d((f_0(a), a),(f_i(b), b))=3$.
Отсюда следует утверждение.

6. \emph{Пусть $A=\{a_1, a_2, a_3, a_4\}$ и $B=\{b_1, b_2, b_3, b_4\}$ --- 
непересекающиеся коклики графа Шрикханде, 
и пусть $\tau$ --- автоморфизм графа Шрикханде, такой, что для любой вершины $s \in B$ 
выполняется: $\tau(s) \in B$, $\tau(s) \ne s$ и  $d(a_i, b_j)=d(a_i, \tau(b_j))$ для любых $i,j=0,1,2,3$. 
Тогда подграф на множестве вершин $A \cup B$ является объединением двух непересекающихся циклов на 4 вершинах.}

Действительно, подграф на множестве вершин $A \cup B$ двудольный, а из леммы \ref{l:koklika} следует, что это регулярный граф степени 2. Тогда это будет либо цикл $C_8$ на 8 вершинах, либо объединением двух непересекающихся циклов $C_4$ на 4 вершинах. 
Но в цикле $C_8$ найдется пара вершин из $A$, имеющая ровно одного общего соседа из $B$, скажем вершину $b$. Так как $\tau$ --- автоморфизм, то $\tau(b)=b$, и мы получаем противоречие.

7. \emph{Пусть $A_0, A_1, A_2, A_3$ --- разбиение графа Шрикханде на 
коклики и для любых $s \ne t$ подграф на множестве вершин $A_s \cup A_t$ является объединением двух непересекающихся циклов на 4 вершинах. 
Тогда $\{A_i: i=0,1,2,3\}$ --- разбиение a) из рис. \ref{f:latin}}.

Докажем, что $L$ --- разбиение a) из рис. \ref{f:latin}. 
Для произвольных $s \ne t$ рассмотрим множества $\{f_0(a): a \in R_s\}$ и $\{f_0(b): b \in R_t\}$.
По пункту 1 эти множества равны $L_s$ и $L_t$ сответственно.
По пукту 4 для любой вершины $a$ из $R_s$ значение $g_0(a)$ равно некоторому $g_1$.
Аналогчино, для любой вершины $b$ из $R_t$ значение $g_0(b)$ равно некоторому $g_2 \ne g_1$.
По пункту 3 для вершины $b \in R_t$ существует единственное $i \ne 0$, такое, что $g_i(b)=g_1$.
Рассмотрим $\tau(s)=f_i(f^{-1}_0(s))$. Из пункта (iv) леммы \ref{l:qwer} следует, что это будет автоморфизм. 
Из пунктов 1, 2 и 5 следует, что к $L_s$, $L_t$ и $\tau$ можно применить пункт 6.
В силу произвольности $s$ и $t$ из пункта 7 следует, что 
$\{L_i: i=0,1,2,3\}$ --- разбиение a) из рис. \ref{f:latin}. 

МДР код с параметрами $(2+2,4^3,4)$ приведен в таблице \ref{tab:2231} Приложения. Обозначим его через $C$.
Также вершины кода обозначим через $(f_i(a),a,i,g_i(a))$.
Докажем, что он единственнен с точностью до эквивалентности.

Пусть $C'$ --- $(2+2,4^3,4)$ МДР кода. Докажем, что он эквивалентен $C$.

Действительно, проекции $C_{;2}$ и $C'_{;2}$ эквивалентны.
Значения $g_0(a)$ для всех $a \in \VV Sh$ по пункту 4 определяются с точностью до перестановки.
Тогда код $C'$ перестановкой и набором автоморфизмов можно перевести в код $C''$, такой, что
если обозначить вершины кода через $(f'_i(a),a,i,g'_i(a))$, то
для любой вершины $a \in \VV Sh$ и любого $i=0,1,2,3$ выполнено $f_i(a)=f'_i(a)$, $g_0(a)=g'_0(a)$.

Докажем, что для любого $i=1,2,3$ и любой вершины $a \in \VV Sh$ выполняется $g_i(a)=g'_i(a)$.
Предположим обратное. Пусть существует вершина $a \in \VV Sh$ и $l \in \{1,2,3\}$ такие, что $g_l(a) \ne g'_l(a)$.
Тогда по пункту 3 существует $j \ne l,0$ такое, что $g'_j(a)=g_l(a)$.
Вершина $a$ принадлежит $R_s$ для некоторго $s=0,1,2,3$.
Из пункта 4 следует, что для некоторого $t$ будет $g_0(c)=g_l(a)$ для всех $c \in R_t$, $t \ne s$.
Тогда по лемме \ref{l:koklika} существует вершина $b \in R_t$, такая, что $d(f_0(a),f_0(b))=1$. 
Так как $g_0(b)=g_l(a)$, то пункту 5 $d(f_0(b),f_l(a))=d(f_0(b),f_0(a))=1$.
Так как $g_0(b)=g'_j(a)$, то по пункту 5 $d(f_0(b),f_j(a))=d(f_0(b),f_0(a))=1$.
Тогда вершина $f_0(b)$ смежна с вершинами $f_0(a), f_l(a), f_j(a)$, то есть вершина $f_0(b)$ смежна с треми вершинами из коклики $L_t$,
и мы получаем противоречие по лемме \ref{l:koklika}.

Тогда код $C''$ совпадает с $C$, а следовательно $C'$ эквивалентен $C$.
\end{proof}

\subsection{МДР коды с параметрами $(3+0, 4^3, 4)$}\label{s:304}

Рассмотрим МДР коды с параметрами $(3+0,4^3,4)$.
Пусть $C$ один из таких кодов.
Для любой вершины $a \in \VV Sh$ грань $C^{a;}_{1;}$ будет $(2+0,4^1,4)$ МДР кодом.
Тогда множество кодовых слов можно представить в виде
$\{(a,f_i(a),g_i(a))
\}$, где $i=0,1,2,3$. 
Также обозначим $F(a)=\{f_i(a):i=0,1,2,3\}$, $G(a)=\{g_i(a):i=0,1,2,3\}$. 
Из кодового расстояния следует, что эти множества будут кокликами  в $Sh$.
\begin{predln}\label{p:shshsh}
Не существует МДР кодов с параметрами $(3+0,4^3,4)$.
\end{predln}
\begin{proof}
Предположим обратное. 
Пусть $C$ --- $(3+0,4^3,4)$ МДР код.
Докажем несколько утверждений.

1. \emph{Расстояние между любыми двумя различными вершинами кода $C$ равно либо 4, либо 6.}

Пусть $a$ и $b$ вершины графа Шрикханде, такие, что $d(a,b)=1$.
Посчитаем сумму $$S=\sum_{i,j=0,1,2,3}d((a,f_i(a),g_i(a)),(b,f_j(b),g_j(b)))=$$
$$16 \cdot d(a,b)+\sum_{i,j=0,1,2,3} d(f_i(a),f_j(b))+\sum_{i,j=0,1,2,3} d(g_i(a),g_j(b))$$
Множества вершин $\{f_i(a):i=0,1,2,3\}$ и $\{f_i(b):i=0,1,2,3\}$ --- 
непересекающиеся коклики. 
По лемме \ref{l:koklika}
каждая вершина из одного множества соединена ровно с двумя вершинами из другого.
Тогда, $$\sum_{i,j=0,1,2,3} d(f_i(a),f_j(b))=4(1+1+2+2)=24.$$
Аналогично, $$\sum_{i,j=0,1,2,3} d(g_i(a),g_j(b))=4(1+1+2+2)=24.$$
Следовательно, $S=64$. 
С другой стороны, так как расстояние 
между двумя различными кодовыми вершинами не меньше 4, то $S \ge 64$, и 
неравенство достигается только если
$d((a,f_i(a),g_i(a)),(b,f_j(b),g_j(b)))=4$ для любых $i,j=0,1,2,3$.

2. \emph{Пусть $a_0,a_1,a_2,a_3$ различные вершины графа Шрикханде такие, что $d(a_2,a_3)=2$, 
а все остальные пары вершин смежны. Тогда множества $F(a_0), F(a_1), F (a_2), F(a_3)$ попарно не пересекаются. Аналогично, попарно не пересекаются множества $G(a_0), G(a_1), G(a_2), G(a_3)$.}

Как следует из кодового расстояния, если вершины $a$ и $b$ смежны, то
множества
$F(a)$ и $F(b)$ не пересекаются.
Поэтому остается доказать, что не пересекаются множества $F(a_2)$ и $F(a_3)$.
Предположим обратное.
Допустим, что для некторых $k$ и $l$ имеет место $f_k(a_2)=f_l(a_3)$.
По лемме \ref{l:koklika} во множестве $F(a_0)$ ровно две вершины, обозначим их
$b_1$ и $b_2$,  такие, что $d(f_k(a_2),b_1)=d(f_k(a_2),b_2)=2$, и во множестве $F(a_1)$
ровно две вершины, обозначим их $b_3$ и $b_4$, такие, что $d(f_k(a_2),b_3)=d(f_k(a_2),b_4)=2$.
Тогда по пункту 1  для каждого $i=0,1,2,3$ выполнено $d(g_k(a_2),b_i)=d(g_l(a_3),b_i)=1$, и вершины
$g_k(a_2)$ и $g_l(a_3)$ имеют не менее 4 общих соседей, а так как в графе Шрикханде
любые две различные вершины имеют ровно 2 общих соседа, то $g_k(a_2)$ и $g_l(a_3)$ совпадают.
Но тогда $d((a_2,f_k(a_2),g_k(a_2)),(a_3,f_l(a_3),g_l(a_3)))=2$, и мы получаем противоречие с кодовым расстоянием.
Утверждение для множеств $G$ доказывется аналогично.

3. \emph{Множества $F(a)$ и $F(b)$, а также  множества $G(a)$ и $G(b)$ совпадают тогда и только тогда, когда  $a$ и $b$ принадлежат одной линейной коклике.}
  
Пусть $a, b, c$ --- три попарно смежные вершины. Тогда по пункту 2 для вершины $d$, смежной с двумя из этих вершин, множество $F(d)$ однозначно восстанавливается по множествам $F(a), F(b), F(c)$. 

Тогда если обозначить
$A$=$F(00)$,
$B$=$F(10)$,
$C$=$F(11)$,
то для любой другой вершины $a$ множество $F(a)$ восстанавливается однозначно, и легко убедиться, что утверждение для множеств $F$ выполняется.

Аналогчино, получается утверждение для множеств $G$.

Также, если поменять порядок координат, то из пункта 3 следует, что 

4. \emph{для любой вершины $a$ коклики $F(a)$ и $G(a)$ --- линейные .}

По пукнту 3 для любой вершины $a$ из $Sh$ множества $(C_{3;})^{a;}_{2;}$ и $(C_{2;})^{a;}_{3;})$ будут линейными кокликами.
Кодовые вершины кода $C$ можно представить также в виде $(f'_i(a), a, g'_i(a))$, где $i=0,1,2,3$, $a \in \VV Sh$.
Тогда, аналогично, можно доказать, что множества $(C_{3;})^{a;}_{1;}$ и $(C_{1;})^{a;}_{3;}$ будут линейными кокликами для всех $a \in \VV Sh$.
Также кодовые слова можно представить в виде $(f''_i(a), g''_i(a), a)$, $i=0,1,2,3$, $a \in \VV Sh$.
Тогда, аналогично, можно доказать, что множества $(C_{1;})^{a;}_{2;}$ и $(C_{2;})^{a;}_{1;}$ будут линейными кокликами для всех $a \in \VV Sh$.
Так как $F(a)=(C_{3;})^{a;}_{1;}$, $G(a)=(C_{2;})^{a;}_{1;}$, то утверждение доказано.


Пусть $a_0, a_1, a_2, a_3$ различные вершины графа Шрикханде, такие что $d(a_2,a_3)=2$, 
а все остальные пары вершин смежны.
По пунктам 2 и  4 коклики $F(a_0)$, $F(a_1)$, $F(a_2)$, $F(a_3)$ --- линейные и попарно не пересекаются.
Тогда множество вершин $\{f_i(a_j):i,j=0,1,2,3\}$ образует граф Шрикханде.
Для некоторых $i,j=0,1,2,3$ вершины $f_i(a_2)$ и $f_j(a_3)$ несмежны.
Эти вершины имеют 2 общих соседа, некоторые $u$ и $v$.
Очевидно, что $u, v$ не принадлежат $F(a_2)$ и $F(a_3)$.
Так как вершины $f_i(a_2)$ и $f_j(a_3)$ принадлежат разным линейным кокликам, то элемент 
$(f_i(a_2)-f_j(a_3))$ не принадлежит множеству $\{02,20,22\}$.
Тогда по лемме \ref{l:s}
вершины $u$ и $v$ смежны. Тогда они не могут принадлежать одной и той же коклике, а значит 
$u=f_k(a_0)$, $v=f_l(a_1)$ для некоторых $k, l$.
Тогда по пункту 1 вершины $g_k(a_0), g_l(a_1), g_i(a_2), g_j(a_3)$ попарно несмежны. 
Также по пунктам 2 и 4 коклики $G(a_0), G(a_1), G_(a_2), G_(a_3)$ --- линейные и не пересекаются.
Также множество вершин $\{g_i(a_j):i,j=0,1,2,3\}$ образуют граф Шрикханде.
Вершины $g_k(a_0), g_l(a_1), g_i(a_2), g_j(a_3)$ образуют коклику в этом графе. 
Но в любой коклике существует такая пара вершин $a$ и $b$, что элемент $(a-b) \in \{02,20,22\}$, то есть 
$a$ и $b$ принадлежат одной линейной коклике, но вершины из множества $\{g_k(a_0), g_l(a_1), g_i(a_2), g_j(a_3)\}$ принадлежат различным линейным кокликам, следовательно, это множество вершин не может быть кокликой. Это противоречие доказывает предложение.
\end{proof}

\section{Параметры при которых МДР кодов не существует}\label{s:notexist}

\begin{predln}\label{p:n6d5}
При $2m+n=6$ не существует $(m+n,4^2,5)$ МДР кодов.
\end{predln}
\begin{proof}
Допустим такой код существует.
Обозначим вершины кода через $(a, f(a), g(a))$, 
где $a \in \VV Sh$, а $f(a), g(a)$ --- либо вершина графа $Sh$, либо вершина графа $K^2$.
Определим на множестве вершин $V=\VV Sh$ графы 
$G_1=(V,E_1)$, $G_2=(V,E_2)$, $G_3=(V,E_3)$, где:

$E_1=\{(a,b): a, b \in V, d(a,b)=1\}$;

$E_2=\{(a,b): a, b \in V, d(f(a),f(b))=1\}$;

$E_3=\{(a,b): a, b \in V, d(g(a), g(b))=1\}$.

Множества ребер $E_1, E_2, E_3$ попарно не пересекаются (иначе получим противоречие с кодовым расстоянием).
Каждый из графов $G_1, G_2, G_3$ является либо графом Шрикханде, либо графом $K^2$.

Тогда $|E_1 \cup E_2 \cup E_3|=144$, но в графе на 16 вершинах не больше 120 ребер, 
и мы получаем противоречие.
\end{proof}

\begin{lemman}\label{l:d=2mn1} 
При $2m+n>5$ не существует МДР кодов с параметрами $(m+n,4^2,2m+n-1)$.
\end{lemman}
\begin{proof} 
Предположим обратное.
Пусть $2m+n>5$ и $C$ ---  $(m+n,4^2,2m+n-1)$ МДР код.
Зафиксируем некоторый набор координат 
$(i_1,\ldots,i_v;j_1,\ldots,j_w)$ 
такой, что $2v+w=2m+n-6$.
Тогда по лемме \ref{l:prgr}  проекция   $C_{i_1,\ldots,i_v;j_1,\ldots,j_w}$ будет
$((m-v)+(n-w),4^2,2(m-v)+(n-w)-1)$ МДР кодом, и так как $2(m-v)+(n-w)=6$, то
применив предложение \ref{p:n6d5}, мы получим противоречие. 
\end{proof}

\begin{lemman}\label{l:grham}
Если $2m+n \ge 6$, то не существует $(m+n, 4^{2m+n-2}, 3)$ МДР кодов.
\end{lemman}
\begin{proof}
Пусть $C$ --- $(m+n, 4^{2m+n-2}, 3)$ МДР код.
Для каждой вершины $a \in \VV D(m,n)$ обозначим 
$B_1(a)=\{b \in \VV D(m,n): d(a,b) \le 1\}$.
Тогда для любой вершины $a$ мощность $|B_1(a)|=3(2m+n)+1$.
Из кодового расстояния следует, что для любых двух различных кодовых вершин $a, b$ множества 
$B_1(a)$ и $B_1(b)$ не пересекаются.
Тогда из неравенства $$4^{2m+n-2}(3(2m+n)+1) \le 4^{2m+n}$$ получаем, что 
$2m+n \le 5$.
\end{proof}.

\begin{lemman}\label{l:d=4}
Если $2m+n>6$, то не существует $(m+n, 4^{2m+n-3}, 4)$ МДР кодов.
\end{lemman}
\begin{proof}
Предположим обратное.
Пусть $2m+n>6$ и
$C$ --- $(m+n,4^{2m+n-3},4)$ МДР код.
Если $n>0$, то проекция $C_{;1}$ --- $(m+(n-1), 4^{2m+n-3}, 3)$ МДР код, и утверждение верно по лемме \ref{l:grham}.
Пусть $n=0$.
По лемме \ref{l:kk} существует кодовая вершина $(s_1,\ldots,s_m)$, такая, что
$s_1=\ldots=s_{m-3}=00$.
Рассмотрим грань $C^{00,\ldots,00}_{1,\ldots,m-3}$.
Это будет $(3+0,4^3,4)$ МДР код, и мы получаем противоречие по предложению \ref{p:shshsh}.
\end{proof}.

\begin{predln}\label{p:notexist}
Если $2m+n>6$ и $2<d<2m+n$, то  не существует $(m+n, 4^k, d)$ МДР кодов.
\end{predln}
\begin{proof}
Предположим обратное.
Пусть $2m+n>6$ и $C$ --- $(m+n,4^k,d)$ МДР код. Если $d=3$ либо $d=4$, то утверждение верно по леммам \ref{l:grham} и 
\ref{l:d=4} соответственно. 

Пусть $d>5$.

Для начала рассмотрим случай, когда $n=0$ и $k$ --- нечетно.
Тогда $d$ будет четным и $d \le 2m-2$. 
Зафиксируем некоторый набор координат 
$(i_1,\ldots,i_v)$, такой 
что $2v=d-4$.
Тогда $C_{i_1,\ldots,i_v;}$ будет $((m-v)+0,4^{k},4)$ МДР кодом 
и $2(m-v) \ge 6$. Если $2(m-v)=6$, то мы получаем противоречие по предложению \ref{p:shshsh}.
Если $2(m-v)>6$, то мы получаем противоречие по лемме \ref{l:d=4}.

Если $n>0$ либо $k$ --- четно, то мы можем зафисировать такой набор координат 
$(i_1,\ldots,i_v;j_1,\ldots,j_w)$, 
что $2v+w=k-2$.
Тогда по лемме \ref{l:prgr} грань $C^{00,\ldots,00;0,\ldots,0}_{i_1,\ldots,i_v;j_1,\ldots,j_w}$ ---
МДР код с параметрами
$((m-v)+(n-w),4^2,d)$. 
Так как $d=2m+n-k+1=2(m-v)+(n-w)-1$, 
то по лемме \ref{l:d=2mn1} и предложению \ref{p:n6d5} выполняется
$2(m-v)+(n-w)$ $\le 5$, и $d \le 4$, и мы получаем противоречие с $d>5$.
\end{proof}

\providecommand\href[2]{#2} \providecommand\url[1]{\href{#1}{#1}}
  \providecommand\bblmay{May} \providecommand\bbloct{October}
  \providecommand\bblsep{September} \def\DOI#1{{\small {DOI}:
  \href{http://dx.doi.org/#1}{#1}}}\def\DOIURL#1#2{{\small{DOI}:
  \href{http://dx.doi.org/#2}{#1}}}\providecommand\bbljun{June}

	
\section{Приложение}\label{s:appendix}

\begin{table}[h]
\begin{center}
\begin{tabular}{|c|c|c|c|c|c|c|c|c|c|c|c|c|c|c|c|c|}
\hline
$a$    & 00 & 02 & 20 & 22 & 01 & 03 & 21 & 23 & 10 & 30 & 12 & 32 & 11 & 13 & 31 & 33\\
\hline
$f(a)$ & 00 & 02 & 20 & 22 & 23 & 21 & 03 & 01 & 12 & 32 & 10 & 30 & 31 & 33 & 11 & 13\\
\hline
\end{tabular}
\end{center}
\caption{Первый из двух неэквивалентных $(2+0,4^2,3)$ МДР кодов. Кодовые слова представлены в виде $(a, f(a))$, где $a, f(a) \in \VV Sh$.}
\label{tab:2021}
\end{table}

\begin{table}[h]
\begin{center}
\begin{tabular}{|c|c|c|c|c|c|c|c|c|c|c|c|c|c|c|c|c|}
\hline
$a$    & 00 & 02 & 20 & 22 & 01 & 03 & 21 & 23 & 10 & 12 & 31 & 33 & 11 & 13 & 30 & 32\\
\hline
$f(a)$ & 00 & 02 & 23 & 21 & 22 & 20 & 01 & 03 & 32 & 30 & 10 & 12 & 13 & 11 & 31 & 33\\
\hline
\end{tabular}
\end{center}
\caption{Второй из двух неэквивалентных $(2+0,4^2,3)$ МДР кодов. Кодовые слова представлены в виде $(a, f(a))$, где $a, f(a) \in \VV Sh$.}
\label{tab:2022}
\end{table}

\begin{table}[h]
\begin{center}
\begin{tabular}{|c|c|c|c|c|c|c|c|c|c|c|c|c|c|c|c|c|}
\hline
$a$    & 00  & 02  & 21  & 23  & 01  & 03  & 22  & 20  & 10  & 12  & 31  & 33  & 11  & 13  & 32  & 30\\
\hline
$f(a)$ & 00  & 12  & 23  & 31  & 21  & 33  & 10  & 02  & 11  & 03  & 30  & 22  & 32  & 20  & 01  & 13\\
\hline
\end{tabular}
\end{center}
\caption{Единственный с точностью до эквивалентности $(1+2,4^2,3)$ МДР код. Кодовые слова представлены в виде $(a, f(a))$, где $a \in \VV Sh$, $f(a) \in \VV K^2$.}
\label{tab:1221}
\end{table}

\begin{table}[h]
\begin{center}
\begin{tabular}{|c|c|c|c|c|c|c|c|c|c|c|c|c|c|c|c|c|}
\hline
$a$    & 00 & 02 & 20 & 22 & 01 & 03 & 21 & 23 & 10 & 30 & 12 & 32 & 11 & 13 & 31 & 33\\
\hline
$f(a)$ & 00 & 02 & 20 & 22 & 23 & 21 & 03 & 01 & 12 & 32 & 10 & 30 & 31 & 33 & 11 & 13\\
\hline
$g(a)$ & 0  & 0  & 0  & 0  & 1  & 1  & 1  & 1  & 2  & 2  & 2  & 2  & 3  & 3  & 3  & 3\\
\hline
\end{tabular}
\end{center}
\caption{Первый из двух неэквивалентны $(2+1,4^2,4)$ МДР кодов. Кодовые слова представлены в виде $(a, f(a), g(a))$, где $a, f(a) \in \VV Sh$, $g(a) \in \VV  K$.}
\label{tab:2121}
\end{table}

\begin{table}[h]
\begin{center}
\begin{tabular}{|c|c|c|c|c|c|c|c|c|c|c|c|c|c|c|c|c|}
\hline
$a$    & 00 & 02 & 20 & 22 & 01 & 03 & 21 & 23 & 10 & 12 & 31 & 33 & 11 & 13 & 30 & 32\\
\hline
$f(a)$ & 00 & 12 & 23 & 31 & 21 & 33 & 10 & 02 & 11 & 03 & 30 & 22 & 32 & 20 & 01 & 13\\
\hline
$g(a)$ & 0  & 0  & 0  & 0  & 1  & 1  & 1  & 1  & 2  & 2  & 2  & 2  & 3  & 3  & 3  & 3\\
\hline
\end{tabular}
\end{center}
\caption{Второй из двух неэквивалентных $(2+1,4^2,4)$ МДР кодов. Кодовые слова представлены в виде $(a, f(a), g(a))$, где $a, f(a) \in \VV Sh$, $g(a) \in K$.}
\label{tab:2122}
\end{table}

\begin{table}[h]
\begin{center}
\begin{tabular}{|c|c|c|c|c|c|c|c|c|c|c|c|c|c|c|c|c|}
\hline
$a$    & 00 & 02 & 21 & 23 & 01 & 03 & 22 & 20 & 10 & 12 & 31 & 33 & 11 & 13 & 32 & 30\\
\hline
$f(a)$ & 00 & 12 & 23 & 31 & 21 & 33 & 10 & 02 & 11 & 03 & 30 & 22 & 32 & 20 & 01 & 13\\
\hline
$g(a)$ & 0  & 0  & 0  & 0  & 1  & 1  & 1  & 1  & 2  & 2  & 2  & 2  & 3  & 3  & 3  & 3\\
\hline
\end{tabular}
\end{center}
\caption{Единственный с точностью до эквивалентности $(1+3,4^2,4)$ МДР кодов. Кодовые слова представлены в виде $(a, f(a))$, где $a \in \VV Sh$, $f(a) \in \VV K^2$, $g(a) \in K$.}
\label{tab:1321}
\end{table}

\begin{table}[h]
\begin{center}
\begin{tabular}{|c|c|c|c|c|c|c|c|c|c|c|c|c|c|c|c|c|}
\hline
$a$      & 00 & 02 & 20 & 22 & 01 & 03 & 21 & 23 & 10 & 30 & 12 & 32 & 11 & 13 & 31 & 33\\
\hline
$i=0$    & 00 & 02 & 20 & 22 & 23 & 21 & 03 & 01 & 12 & 32 & 10 & 30 & 31 & 33 & 11 & 13\\
\hline
$i=1$    & 02 & 00 & 22 & 20 & 21 & 23 & 01 & 03 & 10 & 30 & 12 & 32 & 33 & 31 & 13 & 11\\
\hline
$i=2$    & 20 & 22 & 00 & 02 & 03 & 01 & 23 & 21 & 32 & 12 & 30 & 10 & 11 & 13 & 31 & 33\\
\hline
$i=3$    & 22 & 20 & 02 & 00 & 01 & 03 & 21 & 23 & 30 & 10 & 32 & 12 & 13 & 11 & 33 & 31\\
\hline 
\end{tabular}
\end{center}
\caption{Первый из двух неэквивалентных $(2+1,4^3,3)$ МДР кодов. Кодовые слова представлены в виде $(f_i(a), a, i)$, где $f_i(a), a \in \VV Sh$, $i \in \VV K$.
В таблице представлены значения $f_i(a)$.}
\label{tab:2131}
\end{table}

\begin{table}[h]
\begin{center}
\begin{tabular}{|c|c|c|c|c|c|c|c|c|c|c|c|c|c|c|c|c|}
\hline
$a$      & 00 & 02 & 20 & 22 & 01 & 03 & 21 & 23 & 10 & 30 & 12 & 32 & 11 & 13 & 31 & 33\\
\hline
$i=0$    & 00 & 02 & 23 & 21 & 22 & 20 & 01 & 03 & 32 & 30 & 10 & 12 & 13 & 11 & 31 & 33\\
\hline
$i=1$    & 02 & 00 & 21 & 23 & 20 & 22 & 03 & 01 & 30 & 32 & 12 & 10 & 11 & 13 & 33 & 31\\
\hline
$i=2$    & 21 & 23 & 00 & 02 & 01 & 03 & 20 & 22 & 12 & 10 & 32 & 30 & 33 & 31 & 13 & 11\\
\hline
$i=3$    & 23 & 21 & 02 & 00 & 03 & 01 & 22 & 20 & 10 & 12 & 30 & 32 & 31 & 33 & 11 & 13\\
\hline 
\end{tabular}
\end{center}
\caption{Второй из двух неэквивалентных $(2+1,4^3,3)$ МДР кодов. Кодовые слова представлены в виде $(f_i(a), a, i)$, где $f_i(a), a \in \VV Sh$, $i \in \VV K$. В таблице представлены значения $f_i(a)$.}
\label{tab:2132}
\end{table}

\begin{table}[h]
\begin{center}
\begin{tabular}{|c|c|c|c|c|c|c|c|c|c|c|c|c|c|c|c|c|}
\hline
$a$      & 00  & 12  & 23  & 31  & 21  & 33  & 10  & 02  & 11  & 03  & 30  & 22  & 32  & 20  & 01  & 13\\
\hline   
$i=0$    & 00  & 02  & 21  & 23  & 01  & 03  & 22  & 20  & 10  & 12  & 31  & 33  & 11  & 13  & 32  & 30\\
\hline
$i=1$    & 02  & 00  & 23  & 21  & 03  & 01  & 20 &  22  & 12  & 10  & 33  & 31  & 13  & 11  & 30  & 32\\
\hline
$i=2$    & 21  & 23  & 02  & 00  & 22  & 20  & 03  & 01  & 31  & 33  & 12  & 10  & 32  & 30  & 13  & 11\\
\hline
$i=3$    & 23  & 21  & 00  & 02  & 20  & 22  & 01  & 03  & 33  & 31  & 10  & 12  & 30  & 32  & 11  & 13\\
\hline 
\end{tabular}
\end{center}
\caption{Единственный с точностью до эквивалентности $(2+1,4^3,3)$ МДР код. Кодовые слова представлены в виде $(f_i(a), a, i)$, $f_i(a) \in \VV Sh$, $a \in \VV K^2$, $i \in \VV K$.
В таблице представлены значения $f_i(a)$.
}
\label{tab:1331}
\end{table}

\begin{table}[h]
\begin{center}
\begin{tabular}{|c@{\,}|@{\,}c@{\,}|@{\,}c@{\,}|@{\,}c@{\,}|@{\,}c@{\,}|@{\,}c@{\,}|@{\,}c@{\,}|@{\,}c@{\,}|@{\,}c@{\,}|@{\,}c@{\,}|@{\,}c@{\,}|@{\,}c@{\,}|@{\,}c@{\,}|@{\,}c@{\,}|@{\,}c@{\,}|@{\,}c@{\,}|@{\,}c@{\,}|}
\hline
$a$      & 00   & 02   & 20   & 22   & 01   & 03   & 21   & 23   & 10   & 30   & 12   & 32   & 11   & 13   & 31   & 33 \\ \hline\hline
$i=0$    & 00\,0 & 02\,0 & 20\,0 & 22\,0 & 23\,1 & 21\,1 & 03\,1 & 01\,1 & 12\,2 & 32\,2 & 10\,2 & 30\,2 & 31\,3 & 33\,3 & 11\,3 & 13\,3\\\hline
$i=1$    & 02\,1 & 00\,1 & 22\,1 & 20\,1 & 21\,0 & 23\,0 & 01\,0 & 03\,0 & 10\,3 & 30\,3 & 12\,3 & 32\,3 & 33\,2 & 31\,2 & 13\,2 & 11\,2\\\hline
$i=2$    & 20\,2 & 22\,2 & 00\,2 & 02\,2 & 03\,3 & 01\,3 & 23\,3 & 21\,3 & 32\,0 & 12\,0 & 30\,0 & 10\,0 & 11\,1 & 13\,1 & 31\,1 & 33\,1\\\hline
$i=3$    & 22\,3 & 20\,3 & 02\,3 & 00\,3 & 01\,2 & 03\,2 & 21\,2 & 23\,2 & 30\,1 & 10\,1 & 32\,1 & 12\,1 & 13\,0 & 11\,0 & 33\,0 & 31\,0\\
\hline 
\end{tabular}
\end{center}
\caption{Единственный с точностью до эквивалентности $(2+1,4^3,3)$ МДР код. Кодовые слова представлены в виде $(f_i(a), a, i, g_i(a))$, где $a, f_i(a) \in \VV Sh$, $i, g_i(a) \in \VV K$.
В таблице представлены значения $f_i(a), g_i(a)$.}
\label{tab:2231}
\end{table}

\end{document}